\input amstex
\documentstyle{amsppt}
\magnification 1200
\NoBlackBoxes
\nologo
\topmatter
\title On Quantum Cohomology Rings of Partial Flag Varieties
\endtitle
\author Ionu\c t Ciocan-Fontanine
\endauthor
\affil Institut Mittag-Leffler \endaffil
\date February 9, 1997  \enddate

\address{Institut Mittag-Leffler, Aurav\" agen 17, S-182 62, 
Djursholm, Sweden} \endaddress
\email{ciocan\@ml.kva.se} \endemail

\endtopmatter
\document
\heading { 0. Introduction} \endheading

The main goal of this paper is to give a unified description for the structure
of the small quantum cohomology rings for all homogeneous spaces of 
$SL_n(\Bbb C)$.

The quantum cohomology ring of a smooth projective variety, or, more generally 
of a symplectic manifold, has been introduced by physicists in the study of
topological field theories (\cite{V}, \cite{W}). 
In the past few years, the highly non-trivial task
of giving a rigorous mathematical treatment for the theory of quantum
cohomology has been accomplished, both in the realm of algebraic and 
symplectic geometry. In various degrees of generality, this can be found in \cite{B}, \cite{BM}, \cite{KM}, \cite{Kon}
\cite {LT1}, \cite {LT2}, \cite{MS}, \cite{RT}, as well as the surveys 
\cite{FP} and \cite{T}.
Roughly speaking, the quantum cohomology ring of a variety $X$ is a deformation of the usual cohomology ring, with parameter space given by $H^*(X)$. The
multiplicative structure of quantum cohomology encodes the enumerative
geometry of rational curves on $X$ 

If one restricts the parameter space to $H^{1,1}(X)$, one gets the {\it small quantum cohomology ring} (terminology taken from \cite{FP}). This ring, in the case of partial flag varieties is the object of the present paper. In order to
state our main results, we will first describe briefly the "classical" side of
the story.

Let $P$ be a parabolic subgroup of $SL_n(\Bbb C)$. We will interpret the 
homogeneous space $F:=SL_n(\Bbb C)/P$ as the complex projective variety parametrizing flags of {\it quotients} of $\Bbb C^n$ of given ranks,
say $n_k>\dots >n_1$. 

By a classical result of C. Ehresmann (\cite {E}),
the integral cohomolgy of $F$ can be described geometrically as the free
abelian group generated by the {\it Schubert classes}.
These are the (Poincar\'e duals of) fundamental classes of certain subvarieties $\Omega_w\subset F$, one for each element of
the subset $S:=S(n_1,\dots ,n_k)$ of the symmetric group $S_n$,
consisting of permutations $w$ such that if $w(i)>w(i+1)$, then 
$i\in\{ n_1,\dots ,n_k\} $. 

A description of the multiplicative structure is provided by yet another classical theorem, due to A. Borel (\cite{Bor}), which gives 
a presentation for $H^*(F,\Bbb Z)$. Specifically, let
$\sigma_1^1,\dots ,\sigma_{n_1}^1,
\sigma_1^2, \dots ,\sigma_{n_2-n_1}^2,\dots ,\sigma_1^{k+1},\dots ,
\sigma_{n-n_k}^{k+1}$ be $n$ independent variables. Define $A_n$ to be the block
diagonal matrix $\operatorname{diag}(D_1,D_2,\dots ,D_{k+1})$, where
$$D_j:=\pmatrix \sigma_1^j &\sigma_2^j &\hdots &\sigma_{n_j-n_{j-1}-1}^j&\sigma_{n_j-n_{j-1}}^j\\
-1& 0 &\hdots &0 &0\\ \vdots &\vdots &\hdots &\vdots &\vdots\\ 
0 &0 &\hdots  &-1 &0\endpmatrix .$$ 
Borel's result then states that there is a canonical isomorphism
$$\Bbb Z[\sigma _1^1,\dots ,\sigma_{n_1}^1,
\sigma_1^2, \dots ,\sigma_{n_2-n_1}^2,\dots ,\sigma_1^{k+1},\dots ,
\sigma_{n-n_k}^{k+1}]/(g_1,g_2,\dots ,g_n)\cong H^*(F),\tag 0.1 $$
where the $g_j$'s are the coefficients of the characteristic polynomial of the matrix $A_n$. 

The natural problem arising from the above descriptions is to look for polynomial representatives for the Schubert classes. The first case in which this problem has been solved is when $F$ is a Grassmannian, and goes back to G. Giambelli. 

The general case was obtained independently by I. N. Bernstein, I. M. Gelfand, and S. I. Gelfand (\cite{BGG}), and M. Demazure (\cite{D}). 
In fact, it should be noted that it
suffices to solve the above problem for the {\it complete} flag variety
$F_n=SL_n(\Bbb C)/B$. The point is that the map 
$$H^*(F,\Bbb Z)\longrightarrow H^*(F_n,\Bbb Z)\tag 0.2$$ 
induced by flat pull-back via the natural projection
$F_n\rightarrow F$ is an embedding. To be more precise, the Borel description
for the cohomology of the complete flag variety is
$$\Bbb Z[x_1,x_2,\dots ,x_n]/(e_1,\dots ,e_n)\cong H^*(F_n,\Bbb Z),$$
where $e_j$ is the $j^{th}$ elementary symmetric polynomial in $x_1,\dots ,
x_n$. A particulary nice set of representatives for the Schubert classes
in this case are the {\it Schubert polynomials} $\frak S_w(x_1,\dots ,x_n)$
of Lascoux and Sch\"utzenberger (\cite{LS1}, \cite{LS2}).
If we interpret each $\sigma_i^j$ as the $i^{th}$ elementary
symmetric polynomial in variables $x_{n_{j-1}+1},\dots ,x_{n_j}$, then
the image of $H^*(F,\Bbb Z)$ by the map (0.2) is the subring of polynomials
which are symmetric in variables in each of the groups
$$\underbrace{x_1\dots ,x_{n_1}},\underbrace{x_{n_1+1},\dots ,x_{n_2}},\dots ,
\underbrace{x_{n_k+1},\dots ,x_n}.$$ If $w\in S$, then $\frak S_w$ satisfies the
above symmetry, hence it determines a polynomial $P_w$ in the $\sigma $ variables, which represents $[\Omega_w]$ in $H^*(F,\Bbb Z)$. We will call
the $P_w(\sigma )$'s the {\it Giambelli polynomials} associated to $F$.

Among the Schubert varieties, there are the so-called {\it special Schubert
varieties}, which are geometric realizations of the Chern classes of the
universal quotient bundles on $F$. They correspond to the cyclic permutations
$\alpha_{i,j}:=s_{n_j-i+1}\cdot\dots\cdot s_{n_j}$, $1\leq j\leq k$,
$1\leq i\leq n_j$, where $s_m:=(m,m+1)$ is the simple transposition interchanging $m$ and $m+1$. 
Again, when $F$ is a Grassmannian, there is a classical
formula, due to M. Pieri, expressing the product
$[\Omega_{\alpha_{i,1}}]\cdot [\Omega_w]$ in the basis of Schubert classes.
Its generalization to the case of the complete flag variety, hence, by the
above discussion, to any partial flag variety as well, was first stated by
Lascoux and Sch\"utzenberger (\cite{LS1}), and was given a geometric
proof only recently by F. Sottile (\cite{S}). 

In analogy to the case of Grassmannians, we will refer to the Giambelli and Pieri-type formulae as the {\it classical Schubert Calculus} on $F$.

The small quantum cohomology ring of $F$, denoted by $QH^*(F)$, is defined as the $\Bbb Z[q_1,\dots ,q_k]$-module
$H^*(F,\Bbb Z)\otimes _{\Bbb Z}\Bbb Z[q_1,\dots ,q_k]$, where $q_1,\dots ,q_k$ are formal variables, with a new multiplication, which we denote by $*$,
obtained essentially by replacing the classical structure constants with the
{\it 3-point, genus 0, Gromov-Witten invariants} of $F$. A presentation of
$QH^*(F)$ has been given independently by A. Astashkevich and V. Sadov (\cite{AS}), and B. Kim (\cite{Kim1}), with the proof completed in \cite{Kim2}
(the "extreme" cases of Grassmannians and complete flags were established slightly earlier in \cite{ST}, and \cite{C-F1} and \cite{GK} respectively).
Their result is as follows. Let $B_n=(b_{lm})_{1\leq l,m\leq n}$ be the matrix with entries
$$b_{lm}=\cases (-1)^{n_{j+1}-n_j+1}q_j, & \text{if}\ l=n_{j-1}+1,\ m=n_{j+1},
\ 1\leq j\leq k\\ -1, & \text{if}\ l=n_j+1,\ m=n_j,\ 1\leq j\leq k-1\\ 0, & \text{otherwise}.\endcases$$
Then there exists a canonical isomorphism
$$\Bbb Z[\sigma _1^1,\dots ,\sigma_{n_1}^1,\dots ,\sigma_1^{k+1},\dots ,
\sigma_{n-n_k}^{k+1}][q_1\dots ,q_k]/(G_1,\dots ,G_n)\cong QH^*(F),\tag 0.3 $$
where $G_1,\dots ,G_n$ are the coefficients of the characteristic polynomial of the deformed matrix $A_n^q:=A_n+B_n$. 

From the point of view of enumerative geometry, one is interested in
computing the Gromov-Witten invariants of $F$, and the description (0.3) is
not too helpful, unless one has quantum versions of the Giambelli and
Pieri formulas. In other words, one is interested in developing a
{\it Quantum Schubert Calculus}. The first such formulas, in the
case $F$ is a Grassmannian, were discovered by A. Bertram, whose paper
\cite{Be} pioneered the subject. 
Later on, his approach was extended to the case of complete flags, to
obtain the quantum Giambelli formula for the {\it special} Schubert classes
(see \cite{C-F1}, \cite{C-F2}). Using this, the {\it Quantum Schubert polynomials} were constructed with algebro-combinatorial methods
by S. Fomin, S. Gelfand and A. Postnikov (\cite{FGP}), giving therefore the full quantum Giambelli formula for the variety of complete flags. They have also
given a special case of the Quantum Pieri formula, namely the Quantum Monk
formula, which corresponds to multiplying by the {\it first} Chern class of
one of the tautological bundles.

As opposed to the situation of the classical cohomology, the quantum story
for a partial flag variety is far from being determined by the one for the complete flags, the reason being that the quantum cohomology lacks the functoriality enjoyed by the usual one. The main results of this paper are
unified quantum versions of the Giambelli and Pieri formulas, which hold for
{\it any} $F$. These formulas specialize to the ones described above when
$F$ is either a Grassmannian, or the complete flag variety. In order to state
them, we introduce first some notation.

Let 
$$1\leq h_1<\dots <h_m\leq l_m<\dots <l_1\leq k$$ 
be integers. We denote by
$\bold{h}$, respectively $\bold{l}$, the collections $h_1,\dots ,h_m$
and $l_1,\dots ,l_m$. Let
$$\gamma_{\bold{h}\bold{l}}:=\gamma_{h_m,l_m}\cdot
\gamma_{h_{m-1},l_{m-1}}\cdot\dots\cdot\gamma_{h_1,l_1},$$

$$\delta_{\bold{h}\bold{l}}:=\delta_{h_1,l_1}\cdot\delta_{h_2,l_2}\cdot
\dots\cdot\delta_{h_m,l_m},$$
where $\gamma_{h,l}$ and $\delta_{h,l}$ denote the cyclic permutations $s_{n_h}\cdot\dots\cdot s_{n_{l+1}-1}$ and
$s_{n_l-1}\cdot\dots\cdot s_{n_{h-1}+1}$ respectively,
for any integers $h,l$ satisfying $1\leq h\leq l\leq k$.

Denote by $q_{\bold{h}\bold{l}}$ the monomial 
$$q_{\bold{h}\bold{l}}:=\underbrace{q_{h_1}\dots q_{h_2-1}}\underbrace{q_{h_2}^2\dots q_{h_3-1}^2}\dots 
\underbrace{q_{h_m}^m\dots q_{l_m}^m}\underbrace{q_{l_m-1}^{m-1}\dots q_{l_{m-1}}^{m-1}}\dots\underbrace{ q_{l_2-1}\dots q_{l_1}}.$$
For each $1\leq j\leq k$ and $1\leq i\leq n_j$, let
$\alpha_{i,j}=s_{n_j-i+1}\cdot\dots\cdot s_{n_j}$.
For $1\leq a<b\leq n$ denote by $t_{ab}$ the transposition
interchanging $a$ and $b$.
If $w,w'\in S$, write $w\overset{\alpha_{i,j}}\to{\longrightarrow}w'$ if 
there exist integers $a_1,b_1,\dots ,a_i,b_i$, such that 

(1) $a_r\leq n_j<b_r$, for $1\leq r\leq i$ and $w'=w\cdot t_{a_1b_1}\cdot\dots
\cdot t_{a_ib_i}$;

(2) $\ell (w\cdot t_{a_1b_1}\cdot\dots\cdot t_{a_rb_r})=\ell (w)+r$, 
$1\leq r\leq i$;

(3) the integers $a_1,\dots ,a_i$ are {\it distinct}.

Our first main theorem is the
 
\proclaim{ Quantum Pieri formula } For every $1\leq j\leq k$,
$1\leq i\leq n_j$ and $w\in S$,
$$[\Omega_{\alpha_{i,j}}]*[\Omega_w]=
\sum_{w\overset{\alpha_{i,j}}\to{\longrightarrow}w'}[\Omega_{w'}]+
\sum_{\bold{h},\bold{l}}q_{\bold{h}\bold{l}}
\left ( \sum_{w''}[\Omega_{w''\cdot\delta_{\bold{h}\bold{l}}}]\right ),$$ 
where the second sum is over all collections $\bold{h},\bold{l}$ such that 
$m\leq i$, $h_m\leq j\leq l_m$, and    
$$\ell (w\cdot\gamma_{\bold{h}\bold{l}})=
\ell (w)-\sum_{c=1}^m(n_{l_c+1}-n_{h_c}),$$ 
while the last sum is
over all permutations $w''\in S_n$ satisfying
$w\cdot\gamma_{\bold{h}\bold{l}}\overset{\tilde\alpha_{i,j}}\to{\longrightarrow}w''$, with
$\tilde\alpha_{i,j}=\alpha_{i,j}\cdot s_{n_j}\cdot s_{n_j-1}\cdot\dots\cdot s_{n_j-m+1}$, and 
$$\ell (w''\cdot\delta_{\bold{h}\bold{l}})=\ell (w'')-m-\sum_{c=1}^m(n_{l_c}-n_{h_c-1}). 
\ \ \ \square$$
\endproclaim

For each $1\leq j\leq k$, $1\leq i\leq n_j$, let $g_i^j=g_i^j(\sigma )$ be the polynomial representing the
$i^{th}$ Chern class of the $j^{th}$ universal quotient bundle on $F$.
Alternatively, for each $j$, the polynomials $g_i^j,\ 1\leq i\leq n_j$ 
are the coefficients of the
characteristic polynomial $\operatorname{det}(A_{n_j}+\lambda I)$, where
$A_{n_j}$ is the upper left $n_j\times n_j$ submatrix of the matrix $A_n$.

Define now polynomials $G_i^j=G_i^j(\sigma ,q)$, $1\leq j\leq k\ 1\leq i\leq n_j$  in exactly the same way as above, but using the 
Astashkevich-Sadov-Kim matrix $A_n^q$ instead of $A_n$

For a partition $\Lambda_j:=(\lambda_{j,1},\dots ,\lambda_{j,n_{j+1}-n_j})$ 
with (at most) 
$n_{j+1}-n_j$ parts and such that each part $\lambda_{j,m}$ is at most $n_j$,
set 
$$g_{\Lambda_j}^{(j)}:=g_{\lambda_{j,1}}^jg_{\lambda_{j,2}}^j\dots g_{\lambda_{j,n_{j+1}-n_j}}^j.$$
Define "standard elementary monomials" 
$g_\Lambda:=g_{\Lambda_1\Lambda_2\dots \Lambda_k}\in\Bbb Z[\sigma ]$ by
$$g_\Lambda:=g_{\Lambda_1}^{(1)}g_{\Lambda_2}^{(2)}
\dots g_{\Lambda_k}^{(k)}. $$

Similarly, the {\it standard quantum elementary monomial}
$G_\Lambda$ is the polynomial in $\Bbb Z[\sigma ,q]$ obtained by replacing in $g_\Lambda$ each factor $g_\lambda^j$ by the corresponding $G_\lambda^j$.
It is easy to see that that each Giambelli polynomial can be written uniquely
as a linear combination $P_w=\sum_{\Lambda}a_\Lambda (w) g_\Lambda$, 
with $a_\Lambda (w)$ integers.

Following \cite{FGP}, we define
the {\it quantum Giambelli polynomial} $P_w^q(\sigma ,q)$ by
$$P_w^q(\sigma ,q)=\sum_{\Lambda}a_\Lambda (w) G_\Lambda.$$
We then prove

\proclaim{Quantum Giambelli formula} $[\Omega_w]=P_w^q(\sigma,q)$
in $QH^*(F)$, for all $w\in S$.\endproclaim

We describe now briefly the way our proofs go. For the proof of Quantum Pieri
we use the construction of 3-point Gromov-Witten invariants by means 
of hyperquot schemes. By a certain degeneration technique, the computation of
some of these invariants is reduced to evaluating intersection numbers on 
$F$ itself. All the ideas involved here appear already in the geometric proof
of the quantum Monk formula given in \cite{C-F2}.

Since the above definition of quantum Giambelli polynomials is the straightforward extension of the one given in \cite{FGP} for the complete
flag variety, the very nice and simple proof given there for the quantum Giambelli formula will work in the partial flag case too, once
all the crucial ingredients are in place. More precisely, what needs to be shown is first that the formula holds for the {\it special} Schubert classes, and secondly that the quantum Giambelli polynomials defined above are orthogonal
with respect to a naturally defined inner product on the quotient
$\Bbb Z[\sigma ,q]/I_q$ (see (0.3)). 

The first result follows rather easily from the quantum Pieri formula (cf. also \cite{C-F2}). As a byproduct, we also get an independent proof for the Astashkevich-Sadov-Kim theorem (0.3).

The proof of the orthogonality property given in \cite{FGP} is
combinatorial, and is the technical heart of their paper. Rather than trying
to extend their method to the case of partial flags, we provide here a geometric
proof. Namely, by using the fact that the quantum Giambelli formula holds for the special Schubert classes, we reduce the orthogonality to a statement about
vanishing of certain Gromov-Witten invariants. The later is then shown by a degeneration argument similar to the one in the proof of quantum Pieri.

The paper is divided into two main parts. The first three sections contain
a quick review of the results about the classical and quantum cohomology rings
that we will need, and the proof of the quantum Giambelli formula, assuming that Quantum Pieri and orthogonality hold.
The last three sections study the geometry of hyperquot schemes, from which
we deduce the proofs of Quantum Pieri and orthogonality.

\phantom{X}

\subhead{Acknowledgements} \endsubhead

I have learned the subject from Aaron Bertram, and many of the ideas
that are used in this paper originate in his work on quantum cohomology
of Grassmannians. I am also indebted to William Fulton, and Bumsig Kim, 
for useful discussions during the preparation of the paper. Financial support, 
via a postdoctoral fellowship, as well as a stimulating atmosphere for research, has been provided by the Mittag-Leffler Institute.

\heading{ 1. The Classical Cohomology Ring}
\endheading
\subhead{ 1.1 Schubert varieties}\endsubhead

\phantom{X}

Let $0=n_0<n_1<n_2<\dots <n_k<n_{k+1}=n$ be integers. 
Let $V$ be a complex $n$-dimensional vector space. The data
$k,\ n_j,\ j=0,\dots ,k+1$, and $V$ will be fixed for the rest of the paper.
Define $F:=F(n_1,\dots, n_k,V)$ 
to be the variety parametrizing flags of {\it quotients} of $V$, with ranks
given by the $n_j$'s. 
$F$ is a smooth, irreducible, projective variety, of dimension 
$f:=\sum_{j=1}^k(n-n_j)(n_j-n_{j-1})$.
It comes with a tautological sequence of quotient bundles
$$V_F:=V\otimes\Cal O_F
\twoheadrightarrow Q_k\twoheadrightarrow Q_{k-1}\twoheadrightarrow \dots
\twoheadrightarrow Q_1,$$ with $\operatorname{rank}(Q_j)=n_j$.

Let $S_n$ be the symmetric group on $n$ letters and let
$S:=S(n_1,\dots ,n_k)\subset S_n$ be the subset consisting of permutations
$w$ satisfying the condition: if $w(q)>w(q+1)$, then $q\in\{ n_1,\dots ,n_k\} $.
In other words, when regarded as a function $[1,n]\rightarrow [1,n]$, $w$ is 
increasing on each of the intervals $[1,n_1],[n_1+1,n_2],\dots [n_k,n_{k+1}]$.
The {\it rank function} of a permutation $w\in S_n$ is defined by
$$r_w(q,p)=\operatorname{card}\{ i\mid i\leq q, w(i)\leq p\},\
1\leq q,p\leq n .$$

Fix a complete flag of subspaces $V_1\subset V_2\subset\dots\subset V_{n-1}
\subset V_n=V$. For $w\in S$, the corresponding Schubert variety is defined by
$$\Omega_w:=\{ x\in F \mid \operatorname{rank}_x(V_p\otimes\Cal O
\rightarrow Q_q)\leq r_w(q,p), q\in \{ n_1,\dots ,n_k\} ,1\leq p\leq n\} .$$
$\Omega_w$ is an irreducible subvariety
in $F$, of (complex) codimension equal to the length $\ell (w)$ of
the permutation $w$. 

Throughout the paper $H^*(F)$ will denote the integral cohomolgy of $F$. The following two theorems are classical results 
of C. Ehresmann (\cite{E}, see also \cite{F2, Example 14.7.16}).

\proclaim{Theorem 1.1 (Basis)} 
$\{ [\Omega_w]\} _{w\in S}$  freely generate $H^*(F)$ over $\Bbb Z$. \quad\qed
\endproclaim

\proclaim{Theorem 1.2 (Duality)} 
For every $w\in S$ there exists an unique
permutation $\check w\in S$ such that
$$\int_F[\Omega_w]\cup [\Omega_v]=\cases{1,\ \ \text{if}\ v=\check w}\\
{0,\ \ \text{otherwise}}\endcases . \quad\qed$$
\endproclaim

\phantom{X}

\subhead{ 1.2 A presentation of $H^*(F)$}\endsubhead

\phantom{X}

Consider on $F$ the vector bundles 
$$L_j:=\operatorname{ker}(Q_j\rightarrow Q_{j-1}),$$ 
and let $\sigma_i^j:=c_i(L_j)$, $1\leq i\leq n_j-n_{j-1}$, $1\leq j\leq k+1$.
Let $x_1,\dots ,x_n$ be independent variables. For all $0\leq i\leq m\leq n$,
let $e_i^m$ denote the $i^{th}$ elementary symmetric function in the variables
$x_1\dots ,x_m$. We regard the variables in each of the groups
$$\underbrace{x_1\dots ,x_{n_1}},\underbrace{x_{n_1+1},\dots ,x_{n_2}},\dots ,
\underbrace{x_{n_k+1},\dots ,x_n}$$ as the Chern roots of the bundles $Q_1,
L_2,\dots, L_{k+1}$ respectively. For each $1\leq j\leq k+1$, the 
polynomials $e_i^{n_j}$ can be written as polynomials $g_i^j=g_i^j(\sigma_1^1,
\dots ,\sigma_{n_{k+1}-n_k}^k)$ in the Chern classes of these bundles.
$g_i^j$ has weighted degree $i$, where each $\sigma_m^*$ is assigned degree $m$.
In particular, we have polynomials $g_i^{k+1}$ for $1\leq i\leq n$.

We will denote the polynomial ring $$\Bbb Z[\sigma_1^1,\dots ,\sigma_{n_1}^1,
\sigma_1^2, \dots ,\sigma_{n_2-n_1}^2,\dots ,\sigma_1^{k+1},\dots ,
\sigma_{n-n_k}^{k+1}]$$ by $\Bbb Z[\sigma ]$. With this notation, we can state
another classical result, due to A. Borel (\cite{Bor}).

\proclaim{Theorem 1.3} There is a canonical isomorphism
$\Bbb Z[\sigma ]/I\cong H^*(F)$, where $I$ is the 
ideal generated by $g_1^{k+1},\dots ,g_n^{k+1}$. \quad\qed
\endproclaim

\phantom{X}

\subhead{ 1.3 Classical Schubert Calculus for $F$}\endsubhead

\phantom{X}

Let us recall the Schubert polynomials of Lascoux and Sch\" utzenberger
(\cite{LS1}, \cite{LS2}). Define operators 
$\partial _i,\ i=1,\dots ,n-1$ on $\Bbb Z[x_1,\dots ,x_n]$ by 
$$\partial _iP=\frac {P(x_1,\dots ,x_n)-P(x_1,\dots ,x_{i-1},x_{i+1},x_i,
x_{i+2},\dots ,x_n)}{x_i-x_{i+1}}.$$ 
For any $w\in S_n$, write $w=w_\circ\cdot s_{i_1}\cdot \dots \cdot s_{i_k}$,
with $k=\frac{n(n-1)}{2}-\ell (w)$, where $s_i=(i,i+1)$ is the transposition 
interchanging $i$ and $i+1$, and $w_\circ$ is the permutation of longest
length, given by $w_\circ (j)=n-j+1$, $1\leq j\leq n$. 
The {\it Schubert polynomial} $\frak S _w(x)\in 
\Bbb Z[x_1,\dots ,x_n]$ is defined by 
$$\frak S _w(x)=\partial _{i_k}\circ \dots 
\circ \partial_{i_1}(x_1^{n-1}x_2^{n-2}\dots x_{n-1}).$$ 

It is shown in \cite{M} that 
if $w\in S$, then the corresponding Schubert polynomial is
symmetric in each group of variables
$$\underbrace{x_1\dots ,x_{n_1}},\underbrace{x_{n_1+1},\dots ,x_{n_2}},\dots ,
\underbrace{x_{n_k+1},\dots ,x_n},$$ hence it can be written as a polynomial
$P_w(\sigma )$, of weighted degree $\ell (w)$. We will call these $P_w(\sigma )$
{\it Giambelli polynomials}. The following theorem is due to Bernstein-Gelfand-Gelfand (\cite{BGG}), and Demazure (\cite{D}).

\proclaim{Theorem 1.4 (Giambelli-type formula)}
$[\Omega_w]=P_w(\sigma)$ in $H^*(F)$. \quad\qed\endproclaim 

In particular, consider the cyclic permutations (of length $i$)
$\alpha_{i,j}:=s_{n_j-i+1}\cdot\dots\cdot s_{n_j}$ and
$\beta_{i,j}:=s_{n_j+i-1}\cdot\dots\cdot s_{n_j}$. Note that these 
permutations are in $S$. Their Schubert polynomials
are $\frak S_{\alpha_{i,j}}=e_i^{n_j}$ and  $\frak S_{\beta_{i,j}}=h_i^{n_j}$,
where $h_i^{n_j}$ is the $i^{th}$ {\it complete} symmetric polynomial in
variables $x_1,\dots ,x_{n_j}$. Let $f_i^j(\sigma)$ be the polynomial in the
$\sigma$-variables obtained from $h_i^{n_j}$. By Theorem 1.3, 
$$[\Omega_{\alpha_{i,j}}]=g_i^j\ \ \text{and}\ \ 
[\Omega_{\beta_{i,j}}]=f_i^j\ \ \text{in}\ H^*(F).\tag {1.1}$$

The following Pieri-type formula, due to
A. Lascoux and M. Sch\"utzenberger \cite{LS1}, has been given recently a geometric proof by F. Sottile \cite{S}:

Let $w,w'\in S$. For $1\leq a<b\leq n$ denote by $t_{ab}$ the transposition
interchanging $a$ and $b$.
Write $w\overset{\alpha_{i,j}}\to{\longrightarrow}w'$ if 
there exist integers $a_1,b_1,\dots ,a_i,b_i$, satisfying 

(1) $a_m\leq n_j<b_m$, for $1\leq m\leq i$ and $w'=w\cdot t_{a_1b_1}\cdot\dots
\cdot t_{a_ib_i}$;

(2) $\ell (w\cdot t_{a_1b_1}\cdot\dots\cdot t_{a_mb_m})=\ell (w)+m$, 
$1\leq m\leq i$;

(3$\alpha $) the integers $a_1,\dots ,a_i$ are {\it distinct}.

Similarly, $w\overset{\beta_{i,j}}\to{\longrightarrow}w'$ if
there exist $a_1,b_1,\dots ,a_i,b_i$ as above, satisfying (1), (2) and

(3$\beta $) $b_1,\dots ,b_i$ are {\it distinct}.

\proclaim{Theorem 1.5 (Pieri-type formula)} The following hold in $H^*(F)$:

$(i)\ \ [\Omega_{\alpha_{i,j}}]\cdot[\Omega_w]=
\sum_{w\overset{\alpha_{i,j}}\to{\longrightarrow}w'}[\Omega_{w'}]$.

$(ii)\ \ [\Omega_{\beta_{i,j}}]\cdot[\Omega_w]=
\sum_{w\overset{\beta_{i,j}}\to{\longrightarrow}w'}[\Omega_{w'}]$.
\quad\qed
\endproclaim

\remark{\bf Remark 1.6} From the exact sequence
$$0\longrightarrow L_j\longrightarrow Q_j\longrightarrow Q_{j-1}
\longrightarrow 0,$$
we get 
$$c_i(Q_j)=\sum_{r=0}^{n_j-n_{j-1}}\sigma_r^jc_{i-r}(Q_{j-1}).$$
But one easily sees that $c_i(Q_j)=[\Omega_{\alpha_{i,j}}]$. Using (1.1), it
follows that the polynomials $g_i^j$ satisfy the following recurssion
(which in fact defines them uniquely):
$$g_i^j=\sum_{r=0}^{n_j-n_{j-1}}\sigma_r^jg_{i-r}^{j-1}, \tag {1.2}$$
where, by convention, we set $g_0^{j-1}=1$ and $g_m^{j-1}=0$, if either $m<0$,
or $m>n_{j-1}$.

Also, using the same exact sequence, the relations (1.1), (1.2), 
and the well known identity
$$\left (\sum_{r=0}^{n_{j-1}}e_r^{n_{j-1}}t^r\right )^{-1}=\sum_{p\geq 0}
(-1)^ph_p^{n_{j-1}}t^p,$$
we obtain that the following identity holds in $H^*(F)$:
$$[\Omega_{\alpha_{i,j}}]=\sum_{r=0}^{n_j-n_{j-1}}\left (\sum_{p=0}^r(-1)^p
[\Omega_{\beta_{p,j-1}}]\cdot [\Omega_{\alpha_{r-p,j}}]\right )\cdot
[\Omega_{\alpha_{i-r,j-1}}],\tag{1.3}$$
where, by convention, $[\Omega_{\alpha_{0,m}}]=[\Omega_{\beta_{0,m}}]=1$ and
$[\Omega_{\alpha_{<0,m}}]=0$, for all $m$. \quad\qed
\endremark

\heading {\bf 2. The Small Quantum Cohomology Ring of $F$}\endheading

We give below the precise definition of the small quantum cohomology ring
only for the specific case of a partial flag manifold.

The {\it 3-point, genus 0, Gromov-Witten (GW) invariants} of $F$,
which we denote by
$I_{3,\beta }(\gamma_1\gamma_2\gamma_3)$, are defined as intersection numbers
on Kontsevich's moduli space of stable maps $\overline M_{0,3}(F,\beta )$
(see \cite{KM}, \cite{Kon}, \cite{BM}, \cite{FP}). Here $\beta\in H_2(F)$ and $\gamma_1,\gamma_2,\gamma_3\in H^*(F)$. These numbers have the following enumerative significance (\cite{FP}): 

Let $\Gamma_1,\Gamma_2,\Gamma_3$ be subvarieties of $F$, 
representing the cohomology classes $\gamma_1,\gamma_2,\gamma_3$ respectively.
Let $g_1,g_2,g_3\in SL(n,\Bbb C)$ be general elements, and denote by
$g_i\Gamma_i$ the translate of $\Gamma_i$ by $g_i$. Then 
$I_{3,\beta }(\gamma_1\gamma_2\gamma_3)$ is the number of maps $\mu:\Bbb P^1
\rightarrow F$ such that $\mu_*[\Bbb P^1]=\beta$ and $\mu(\Bbb P^1)$ meets
$g_1\Gamma_1, g_2\Gamma_2$ and $g_3\Gamma_3$. 

Since we will give a different
construction of these invariants in Section 4, we will not say more about them
here. The multiplication in the (small) quantum cohomology ring is defined
using these $I_{3,\beta }$ as structure constants. More precisely, this goes as 
follows:

Introduce formal variables $q_1,\dots ,q_k$, corresponding respectively to
the generators (cf. Theorem 1.1) 
$[\Omega_{\check s_{n_1}}],\dots ,[\Omega_{\check s_{n_k}}]$ of $H_2(F)$.
For a (holomorphic) map $\mu :\Bbb P^1\rightarrow F$ we can write 
$\beta =\mu_*[\Bbb P^1]=\sum_{j=1}^kd_j[\Omega_{\check s_{n_j}}]$, with $d_j$
nonnegative integers. We will say that $\mu$ has multidegree $\overline d=
(d_1,\dots ,d_k)$, and we'll replace $\beta$ by $\overline d$ in the notation
for GW invariants.

Let $K:=\Bbb Z[q_1\dots ,q_k]$. On the $K$-module
$H^*(F)\otimes _{\Bbb Z}K$, define the quantum multiplication $*$ by
putting first
$$[\Omega_u]*[\Omega_v]:=\sum_{\overline d}q_1^{d_1}\dots q_k^{d_k}
\sum_{w\in S}I_{3, \overline d}([\Omega_u][\Omega_v][\Omega_w])
[\Omega_{\check w}],\tag{2.1}$$ for all $u,v\in S$, and then extending 
linearly on $H^*(F)$ and trivially on $K$. The following theorem is a
particular case of the general results on associativity of quantum 
cohomology (\cite{B}, \cite{BM}, \cite{KM}, \cite {LT1}, \cite {LT2}, \cite{MS}, \cite{RT}).

\proclaim{Theorem 2.1} The operation $*$ defines an associative and commutative
$K$-algebra structure on $H^*(F)\otimes _{\Bbb Z}K$. \quad\qed
\endproclaim

$H^*(F)\otimes _{\Bbb Z}K$ together with this multiplication is called the
small quantum cohomology ring of $F$, and denoted by $QH^*(F)$.
The goal of this paper is to give a description analogous to that in Section 1
for this new algebra.

\heading { 3. Quantum Schubert Calculus}\endheading

\subhead{3.1 The quantum version of the Pieri-type formula}\endsubhead

\phantom{X}

We introduce first some notation.
For integers  $h,l$ satisfying $1\leq h\leq l\leq k$, consider the cyclic permutations $\gamma_{h,l}:=s_{n_h}\cdot\dots\cdot s_{n_{l+1}-1}$ and
$\delta_{h,l}:=s_{n_l-1}\cdot\dots\cdot s_{n_{h-1}+1}$. 
Let now $1\leq j\leq k$ and $1\leq i\leq n_j$ be fixed, and let $$m\leq i,\ \
1\leq h_1<\dots <h_m\leq j\leq l_m<\dots <l_1\leq k$$ be integers. 
We denote by
$\bold{h}$ and $\bold{l}$, respectively the collections $h_1,\dots ,h_m$
and $l_1,\dots ,l_m$. Let
$$\gamma_{\bold{h}\bold{l}}:=\gamma_{h_m,l_m}\cdot
\gamma_{h_{m-1},l_{m-1}}\cdot\dots\cdot\gamma_{h_1,l_1},$$

$$\delta_{\bold{h}\bold{l}}:=\delta_{h_1,l_1}\cdot\delta_{h_2,l_2}\cdot
\dots\cdot\delta_{h_m,l_m}.$$
Denote by $q_{\bold{h}\bold{l}}$ the monomial 
$$q_{\bold{h}\bold{l}}:=q_{h_1}\dots q_{h_2-1}q_{h_2}^2\dots
q_{h_3-1}^2\dots q_{h_m}^m\dots q_{l_m}^mq_{l_m-1}^{m-1}\dots q_{l_{m-1}}^{m-1}\dots q_{l_2-1}\dots q_{l_1}.$$
 
\proclaim{Theorem 3.1 (Quantum Pieri formula)} For every $1\leq j\leq k$,
$1\leq i\leq n_j$ and $w\in S$,
$$[\Omega_{\alpha_{i,j}}]*[\Omega_w]=
\sum_{w\overset{\alpha_{i,j}}\to{\longrightarrow}w'}[\Omega_{w'}]+
\sum_{\bold{h},\bold{l}}q_{\bold{h}\bold{l}}
\left ( \sum_{w''}[\Omega_{w''\cdot\delta_{\bold{h}\bold{l}}}]\right ),
\tag {3.1}$$ 
where the second sum is over all collections $\bold{h},\bold{l}$ such that 
$$\ell (w\cdot\gamma_{\bold{h}\bold{l}})=
\ell (w)-\sum_{c=1}^m(n_{l_c+1}-n_{h_c}),$$ 
while the last sum is
over all permutations $w''\in S_n$ satisfying
$w\cdot\gamma_{\bold{h}\bold{l}}\overset{\tilde\alpha_{i,j}}\to{\longrightarrow}w''$, with
$\tilde\alpha_{i,j}=\alpha_{i,j}\cdot s_{n_j}\cdot s_{n_j-1}\cdot\dots\cdot s_{n_j-m+1}$, and 
$$\ell (w''\cdot\delta_{\bold{h}\bold{l}})=\ell (w'')-m-\sum_{c=1}^m(n_{l_c}-n_{h_c-1}). 
$$
\endproclaim

\remark{\bf Remark 3.2} $(i)$ The first term in the right-hand side of the 
formula is the "classical" one, given by Theorem 1.5.

$(ii)$ The condition
$$\ell (w\cdot\gamma_{\bold{h}\bold{l}})=
\ell (w)-\sum_{c=1}^m(n_{l_c+1}-n_{h_c}),$$
can be rephrased equivalently as 
$$w(n_{h_c})>\operatorname{max}\{ w(n_{h_c}+1),\dots ,w(n_{l_c+1})\} ,\ \ 
1\leq c\leq m.$$

$(iii)$ Note that if $m<i$, then
$\tilde\alpha_{i,j}=s_{n_j-i+1}\cdot\dots\cdot s_{n_j-m}$ 
gives the same kind of cyclic permutation
as $\alpha_{i,j}$, but it determines a Schubert variety only on the flag
varieties for which one of the quotients has rank $n_j-m$ !! In fact, as it 
will be seen in the proof of the Theorem, the last sum comes from applying Theorem 1.5 on a flag variety as above. It can be
easily checked however, that the permutations 
$w''\cdot\delta_{\bold{h}\bold{l}}$ are in fact
in $S$, i.e. , they define Schubert varieties on our original 
$F(n_1,\dots,n_k,V)$. (If $m=i$, then $\tilde\alpha_{i,j}$ is the identity
permutation.) 
Also note that for the terms appearing in the last sum we have
$\ell (w''\cdot\delta_{\bold{h}\bold{l}})=
\ell (w)+i-\sum_{c=1}^m(n_{l_c}-n_{h_c-1})-\sum_{c=1}^m(n_{l_c+1}-n_{h_c})$. 

$(iv)$ In the case when $F$ is the {\it complete} flag variety, a Quantum Pieri
formula is stated in the recent preprint \cite{KiMa} of Kirillov and Maeno, and an algebraic proof is suggested. Their formulation is quite different, and we
have not checked if it agrees with what Theorem 3.1 above says in that case.
\quad\qed
\endremark

\phantom{X}

We will prove the above Theorem in Section 6. For the moment, let us see what it says in some special cases.

\phantom{X}

\flushpar
$\bullet$ {\it Grassmannians:} Let $k=1$, $n_1=m$, i.e., $F=G(m,n)$, the
Grassmannian of $m$-dimensional quotients of $V$. Let $w$ be a Grassmannian
permutation of descent $m$ and shape $\lambda=(\lambda_1,\dots ,
\lambda_{m})$,
with $n-m\geq \lambda_1\geq\lambda_2\geq\dots\geq\lambda_{m}\geq 0$. The
partition $\lambda$ is defined by $\lambda_{m-j+1}=w(j)-j$.
Denote $\Omega_{\lambda}:=\Omega_w$. In particular, the subvariety $\Omega_
{\alpha_{i,1}}$ is $\Omega_{(1^i,0^{m-i})}$ with the new notation.
The following result is due to A. Bertram (\cite {Be})

\proclaim{Corollary 3.3 (Quantum Pieri for Grassmannians)} The following  holds in $QH^*(G(m,n))$:
$$[\Omega_{(1^i,0^{m-i})}]*[\Omega_{\lambda}]=(classical\ term)+
q\left (\sum_{\mu }[\Omega_{\mu }]\right ) ,$$ 
where $\mu $ ranges over partitions with
at most $m$ parts, satisfying $|\mu |=|\lambda |+i-n$, and
$\lambda_1-1\geq\mu_1\geq\lambda_2-1\geq\dots\geq\lambda_{m}-1\geq\mu_{m}
\geq 0$.\quad\qed
\endproclaim

\phantom{X}

\flushpar
$\bullet$ {\it Complete flag varieties:} Let $k=n-1$, hence $n_j=j$ for 
all $j$, i.e., $F=F(1,2,\dots ,n-1,V)$. In the case $i=1$, we have 
$\alpha_{1,j}=s_j$ and (3.1) specializes to the Quantum Monk formula of
\cite{FGP} (see also \cite{C-F2}, \cite{Pe}):

\proclaim{Corollary 3.4 (Quantum Monk formula)} One has in $QH^*(F)$
$$[\Omega_{s_j}]*[\Omega_w]=(classical\ term)+\sum_{t_{hl}}q_h\dots q_{l-1}
[\Omega_{w\cdot t_{hl}}],$$ where the sum is over all transposition of
integers $h,l$, with $1\leq h\leq j<l\leq n$, such that $\ell (w\cdot t_{hl})
=\ell (w) -2(l-h)+1$. \quad\qed
\endproclaim

Note, however, that even for the case of complete flags, Theorem 3.1 says
much more than Corollary 3.4 !

Finally, we look now closer to a special case, which will be needed later.
Recall the identity (1.3), which holds in the classical cohomology ring of
our partial flag variety:
$$[\Omega_{\alpha_{i,j}}]=\sum_{r=0}^{n_j-n_{j-1}}\left (\sum_{p=0}^r(-1)^p
[\Omega_{\beta_{p,j-1}}]\cdot [\Omega_{\alpha_{r-p,j}}]\right )\cdot
[\Omega_{\alpha_{i-r,j-1}}].\tag{1.3}$$
We want to compute the right-hand side when the classical product is replaced 
by the quantum product. Of course, the answer is obtained by applying
Theorem 3.1 twice, but this would seem to give, besides the classical term
$[\Omega_{\alpha_{i,j}}]$, lots of "quantum correction" terms. In fact, a
more careful analysis will show that there is either no correction term, or
only one such term which we identify explicitely. It is better to break the
computation into two pieces.

\proclaim{Lemma 3.5} $(i)$ In the classical cohomology ring $H^*(F)$, we have
for $0\leq p\leq r\leq n_j-n_{j-1}$
$$[\Omega_{\beta_{p,j-1}}]\cdot [\Omega_{\alpha_{r-p,j}}]=
[\Omega_{\beta_{p,j-1}\cdot\alpha_{r-p,j}}].$$
$(ii)$ In the quantum cohomology ring $QH^*(F)$ we have
$$[\Omega_{\beta_{p,j-1}}]*[\Omega_{\alpha_{r-p,j}}]=
[\Omega_{\beta_{p,j-1}\cdot\alpha_{r-p,j}}]$$
as well, i.e., there are no quantum correction terms.\endproclaim

\demo {Proof} $(i)$ is a straightforward computation, e.g., using
Sottile's Theorem 1.5. 

$(ii)$ Pick $1\leq h_1<\dots <h_m\leq j\leq l_m<\dots <l_1\leq k$. 
Since $n_j\geq r+n_{j-1}\geq p+n_{j-1}$, we have
also $n_{l_1+1}\geq n_{j+1}\geq n_j+1\geq n_{j-1}+p+1$. Therefore $\beta_{p,j-1}
(n_{l_1+1})>\beta_{p,j-1}(m)$ for all $m<n_{l_1+1}$, 
by the definition of $\beta_{p,j-1}$. In particular
$$\beta_{p,j-1}(n_{l_1+1})>\beta_{p,j-1}(n_{h_1}).\tag 3.2 $$
To get a quantum contribution for the chosen $h_i$ and $l_i$, 
we should have necessarily, by Remark 3.2 $(ii)$,
$$\beta_{p,j-1}(n_{h_1})>\operatorname{max}
\{ \beta_{p,j-1}(n_{h_1}+1),\dots ,\beta_{p,j-1}(n_{l_1+1})\} .$$ 
This contradicts (3.2). \quad\qed
\enddemo

\proclaim{Lemma 3.6} The product $[\Omega_{\beta_{p,j-1}\cdot\alpha_{r-p,j}}]*
[\Omega_{\alpha_{i-r,j-1}}]$ has no quantum correction terms, unless
$r=p=n_j-n_{j-1}$ and $i\geq n_j-n_{j-2}$, in which case there is exactly one
such term, namely $q_{j-1}[\Omega_{\alpha_{i-(n_j-n_{j-2}),j-2}}]$.\endproclaim

\demo{ Proof} This time we need to pick 
$1\leq h_1<\dots <h_m\leq j-1\leq l_m<\dots <l_1\leq k$. 
If any of the $h_i$, or $l_i$ are different from $j-1$, the condition 
$$w(n_{h_i})>\operatorname{max}\{ w(n_{h_i}+1),\dots ,w(n_{l_i+1})\} $$ 
of Remark 3.2 $(ii)$
is easily seen to be contradicted for $w:=\beta_{p,j-1}\cdot\alpha_{r-p,j}$.
Hence $m=1$, $h_1=l_1=j-1$, and we need 
$$w(n_{j-1})>\operatorname{max}\{ w(n_{j-1}+1),\dots ,w(n_j)\} .$$
This happens iff $r=p=n_j-n_{j-1}$, i.e., the only case that may give
quantum contributions is the product
$$[\Omega_{\beta_{n_j-n_{j-1},j-1}}]*[\Omega_{\alpha_{i-(n_j-n_{j-1}),j-1}}],$$
for $h=l=j-1$. In this case $\beta_{n_j-n_{j-1},j-1}\cdot\gamma_{h,l}=id$ (the
identity permutation), and the Quantum Pieri formula (3.1) specializes to give the Lemma. \quad\qed
\enddemo

From the identity (1.3) and the two previous Lemmas we get immediately

\proclaim{Corollary 3.7} The following identity holds in $QH^*(F)$:
$$\align &\sum_{r=0}^{n_j-n_{j-1}}\left (\sum_{p=0}^r(-1)^p
[\Omega_{\beta_{p,j-1}}]*[\Omega_{\alpha_{r-p,j}}]\right )*
[\Omega_{\alpha_{i-r,j-1}}]=\\ &[\Omega_{\alpha_{i,j}}]+(-1)^{n_j-n_{j-1}}
q_{j-1}[\Omega_{\alpha_{i-(n_j-n_{j-2}),j-2}}].\quad\qed\endalign$$
\endproclaim

\phantom{X}

\subhead{3.2 The Quantum Giambelli formula}\endsubhead

\phantom{X}

\definition{ Definition 3.8} For $0\leq j\leq k+1,\ 1\leq i\leq n_j$,	
let $G_i^j\in \Bbb Z[\sigma ,q]$ be the polynomials defined in one of the
following equivalent ways.

(3.8.1) Set $G_0^j:=1,\ G_1^j:=g_1^j$, for all $j$. Then $G_i^j$ for $i\geq 2$,
and all $j$ is defined recursively by
$$G_i^j:=(-1)^{n_j-n_{j-1}+1}q_{j-1}G_{i-(n_j-n_{j-2})}^{j-2}+
\sum_{r=0}^{n_j-n_{j-1}}\sigma_r^jG_{i-r}^{j-1}.$$

(3.8.2) For each $1\leq j\leq k+1$, construct a graph as follows:

$\bullet$ choose $j$ vertices and label them $\bold{v}_1,\dots ,\bold{v}_j$;

$\bullet$ for every $1\leq l\leq j-1$, join the vertices $\bold{v}_l$ and
$\bold{v}_{l+1}$ by an edge and give it the label $(-1)^{n_{l+1}-n_l+1}q_l$;

$\bullet$ for every $1\leq l\leq j$, attach $n_l-n_{l-1}$ tails to the vertex
$\bold{v}_l$, with labels $\sigma_1^l,\dots ,\sigma_{n_l-n_{l-1}}^l$ respectively.

Now define $G_i^j$ to be the sum of all monomials obtained by choosing
edges in this graph and forming the product of their labels, such that the total degree of the monomial is $i$, where 
$\operatorname{deg}(q_l)=n_{l+1}-n_{l-1}$ and 
$\operatorname{deg}(\sigma_m^l)=m$, for every $l,m$, and no two of the chosen 
edges share a common vertex. This description has been shown to us by W. Fulton.

(3.8.3) For each $j$, 
the polynomials $G_i^j,\ 1\leq i\leq n_j$ are the coefficients of the
characteristic polynomial $\operatorname{det}(A_{n_j}^q+\lambda I)$, where
$A_{n_j}^q$ is the upper left $n_j\times n_j$ submatrix of the 
Astashkevich-Sadov-Kim matrix $A_n^q$ (\cite{AS}, \cite{Kim1}, \cite{Kim 2}). 
\quad\qed
\enddefinition

It is immediate from any of these descriptions that $G_i^j(\sigma, 0)=
g_i^j(\sigma )$. We are now ready to formulate a special case of the Quantum
Giambelli formula.

\proclaim{Theorem 3.9} $(i)$ $[\Omega_{\alpha_{i,j}}]=G_i^j(\sigma ,q)$ in $QH^*(F)$, for all $1\leq i\leq n_j,\ 0\leq j\leq k$.

$(ii)$ $G_i^{k+1}(\sigma ,q)=0$ in $QH^*(F)$, for all $1\leq i\leq n$.
\endproclaim

\demo{ Proof} Induction on $j$, using the classical Giambelli formula,
the recursion (1.2) satisfied
by the $g_i^j$'s, the identity (1.3), Corollary 3.7, and the recursion (3.8.1)
satisfied by the $G_i^j$'s (cf. the proof of Theorem 5.6 $(i)$ in \cite {C-F2} for the case of complete flags). \quad\qed
\enddemo

\proclaim{Corollary 3.10 (\cite{AS}, \cite{Kim1}, \cite{Kim2})}There is a
canonical isomorphism
$$QH^*(F)\cong\Bbb Z[\sigma ,q]/I_q,$$ where $I_q$ is the ideal $(G_1^{k+1},
\dots ,G_n^{k+1})$.
\endproclaim

\demo{ Proof} Follows from Theorem 3.9 $(ii)$ and \cite{ST, Theorem 2.2}.
\quad\qed\enddemo

\remark{\bf Remark 3.11} Theorem 3.9 $(ii)$ and Corollary 3.10 were formulated independently by A. Astashkevich and V. Sadov (\cite{AS}), and
B. Kim (\cite{Kim1}), with the proof completed in \cite{Kim2}. 
As far as I know, Theorem 3.9 $(i)$
is new here. For the case of complete flags, Theorem 3.9 and Corollary 3.10
were proved first in \cite{C-F1}. \quad\qed\endremark

Following \cite{FGP}, Theorem 3.9, together with a quantum analogue of 
Theorem 1.2 are sufficient to obtain the general Quantum Giambelli formula.
More precisely, this goes as follows:

For a partition $\Lambda_j:=(\lambda_{j,1},\dots ,\lambda_{j,n_{j+1}-n_j})$ 
with (at most) 
$n_{j+1}-n_j$ parts and such that each part $\lambda_{j,m}$ is at most $n_j$,
set 
$$g_{\Lambda_j}^{(j)}:=g_{\lambda_{j,1}}^jg_{\lambda_{j,2}}^j\dots g_{\lambda_{j,n_{j+1}-n_j}}^j.$$
Define "standard elementary monomials" 
$g_\Lambda:=g_{\Lambda_1\Lambda_2\dots \Lambda_k}\in\Bbb Z[\sigma ]$ by
$$g_\Lambda:=g_{\Lambda_1}^{(1)}g_{\Lambda_2}^{(2)}
\dots g_{\Lambda_k}^{(k)}. \tag 3.3 $$
The number of such monomials is 
$$\sharp\{ g_\Lambda\} =\prod_{j=0}^k\pmatrix n_{j+1} \\ n_j\endpmatrix ,$$
which coincides with the rank of $H^*(F)$.
It follows by realizing $F$ as a succesion of Grassmann bundles that the
monomials $\{ g_\Lambda\}$ generate $H^*(F)$ over $\Bbb Z$. Summarizing, we
have the following proposition.

\proclaim{Proposition 3.12} The standard elementary
monomials $\{ g_\Lambda\} $ 
form a linear basis in $H^*(F)$. \quad\qed\endproclaim

Since the Giambelli polynomials $\{ P_w(\sigma )\} _{w\in S}$ 
also form a basis in $H^*(F)$, we can write uniquely
$$P_w=\sum_{\Lambda}a_\Lambda g_\Lambda ,\tag 3.4 $$ 
with $a_\Lambda$ integers (depending, of course, on $w$).
The following definition was given in the case of complete flags in \cite{FGP}.

\definition{ Definition 3.13} The {\it standard quantum elementary monomial}
$G_\Lambda$ is the polynomial in $\Bbb Z[\sigma ,q]$ obtained by replacing in $g_\Lambda$ each factor $g_\lambda^j$ by the corresponding $G_\lambda^j$, 
defined in 3.8.

The {\it quantum Giambelli polynomial} $P_w^q\in \Bbb Z[\sigma ,q]$ is defined
by
$$P_w^q:=\sum_{\Lambda }a_\Lambda G_\Lambda ,\tag 3.5 $$
with $a_\Lambda$ the integers from (3.4). \quad\qed
\enddefinition

The following is immediate from the definitions 3.8 and 3.13.

\proclaim{Proposition 3.14} $(i)$ $P_w^q(\sigma ,q)$ is a weighted homogeneous
polynomial of weighted degree $\ell (w)$, where $\sigma_i^j$ has degree $i$
and $q_j$ has degree $n_{j+1}-n_{j-1}$, for all $1\leq j\leq k$.

$(ii)$ $P_w^q(\sigma ,0)=P_w(\sigma )$

$(iii)$ (cf. \cite{FGP, 3.6-3.7}) $\{ G_\Lambda\} $ and 
$\{ P_w^q\} $ are
linear bases for $QH^*(F)$. \quad\qed\endproclaim

Let $w^{\circ }\in S$ be the longest element, given by
$w(i)=n-n_j+i-n_{j-1}$, for all $n_{j-1}+1\leq i\leq n_j$, $1\leq j\leq k+1$.
Its length is $\ell (w^{\circ })=\sum_{j=1}^k(n-n_j)(n_j-n_{j-1})=
\operatorname{dim}F$, and $[\Omega_{w^{\circ }}]$ is the class of a point in
$H_0(F)$. By the classical Giambelli formula, 
$$[\Omega_{w^{\circ }}]=P_{w^{\circ}}(\sigma )=(\sigma_{n_1}^1)^{n-n_1}
(\sigma_{n_2-n_1}^2)^{n-n_2}\dots (\sigma_{n_k-n_{k-1}})^{n-n_k},$$
while expressing $[\Omega_{w^{\circ }}]$ in the basis $\{ g_\Lambda\} $ yields
$$[\Omega_{w^{\circ }}]=\underset{n_2-n_1\ \text{factors}}\to
{\underbrace{g_{n_1}^1g_{n_1}^1\dots g_{n_1}^1}}\underset{n_3-n_2\ \text{factors}}\to
{\underbrace{g_{n_2}^2g_{n_2}^2\dots g_{n_2}^2}}\dots 
\underset{n_{k+1}-n_k\ \text{factors}}\to
{\underbrace{g_{n_k}^kg_{n_k}^k\dots g_{n_k}^k}}=g_{\Lambda^{\circ}},$$
where 
$$\Lambda^{\circ }=(\Lambda_1,\dots ,\Lambda_k),\ \Lambda_j=
(\underset{n_{j+1}-n_j\ \text{terms}} \to{\underbrace{n_j,\dots ,n_j}}).$$

For a polynomial $R\in\Bbb Z[\sigma ]$ consider the expansion of its coset 
$R(mod\ I)\in \Bbb Z[\sigma ]/I$ in the basis $\{ P_w\} $, and define
$$\langle R\rangle :=\text{coefficient of}\ P_{w^{\circ }}.$$
Alternately, we can expand $R(mod\ I)$ in the basis $\{ g_\Lambda\} $ and take the coefficient of $g_{\Lambda^{\circ }}$. By the classical Giambelli formula
(Theorem 1.3), we can reformulate Theorem 1.2 as

\proclaim{Proposition 3.15} The polynomials $\{ P_w\} $ satisfy the following
orthogonality property:
$$\langle P_wP_v\rangle =\cases{1,\ \ \text{if}\ v=\check w}\\
{0,\ \ \text{otherwise}}\endcases , $$ 
where $\check w$ is the permutation giving the Schubert variety
dual to $[\Omega_w]$.\quad\qed
\endproclaim

Similarly, for $R(\sigma ,q)\in \Bbb Z[\sigma ,q]$ consider the expansion of
$R(mod\ I_q)$ in the basis $\{ P_w^q\} $ (or $\{ G_\Lambda\} $ respectively), 
and define
$$\langle\langle R\rangle\rangle :=\text{coefficient of}\ P_{w^{\circ }}^q
\ (\text{or}\ G_{\Lambda^{\circ}}\ \text{respectively}).$$

\proclaim{Theorem 3.16 (Orthogonality of the quantum Giambelli polynomials)}
$$\langle\langle P_w^qP_v^q\rangle\rangle =\cases{1,\ \ \text{if}\ v=\check w}\\
{0,\ \ \text{otherwise}}\endcases .$$\endproclaim

\demo{ Proof} The proof will be given in Section 6.\quad\qed
\enddemo

\remark{\bf Remark 3.17} For the special case of complete flags, 
the above Theorem is due to \cite{FGP}, and was proved using combinatorial techniques. The proof we will give in this paper is geometric.\quad\qed
\endremark

\proclaim{Theorem 3.18 (Quantum Giambelli formula)} $[\Omega_w]=P_w^q(\sigma,q)$
in $QH^*(F)$, for all $w\in S$.\endproclaim

\demo{ Proof} Having established the special case of Quantum Giambelli
(Theorem 3.9 $(i)$), and the orthogonality of the $P_w^q$'s (Theorem 3.16),
the proof of the Main Theorem in \cite{FGP, Section 4} applies without changes
in our more general case.\quad\qed
\enddemo

\proclaim{Corollary 3.19 (A. Bertram, \cite{Be})} For the Grassmannian $G(m,n)$,
the classical and quantum Giambelli formulae are the same.
\endproclaim
\demo{ Proof} By the definition 3.8, we have $G_i^1=g_i^1$, for all
$1\leq i\leq m$, hence the quantum Giambelli polynomials coincide with the
classical ones.\quad\qed
\enddemo

We have completed the description of $QH^*(F)$, modulo the proofs of Theorems
3.1 and 3.16. The last three sections of the paper 
are devoted to these proofs, which are based on the geometry of
compactifications of spaces of maps $\Bbb P^1\rightarrow F$ given by
hyperquot schemes. Most of the arguments in \cite{C-F2}, where the case of 
complete flags is treated, require little or no changes, therefore we will
refer to the corresponding results in \cite{C-F2} when appropriate, and give
details only as needed.

\heading{ 4. GW-invariants via hyperquot schemes}\endheading

We recall in this section the construction of 3-point, genus 0 GW-invariants
by means of hyperquot schemes. Details can be found in \cite{Be}, \cite{C-F2}.

\phantom{X}

\subhead{ 4.1 $Hom$ and hyperquot schemes}\endsubhead

\phantom{X}

For fixed $\overline d=(d_1,d_2,\dots ,d_k)$, let
$$H_{\overline d}:=Hom_{\overline d}(\Bbb P^1,F)$$ be the moduli space of
holomorphic maps $\mu:\Bbb P^1\rightarrow F$  of multidegree $\overline d$, i.e., such that $\mu_*[\Bbb P^1]=\sum_{j=1}^k
d_j[\Omega_{\check s_{n_j}}]$. 

Since $F$ is a homogeneous space, standard deformation theory
shows that $H_{\overline d}$ is a smooth quasiprojective 
variety of dimension $$h^0(\Bbb P^1,\mu^*T_F)=
\text {dim} F-\mu_*[\Bbb P^1]\cdot (K_F)=
\sum_{j=1}^k(n-n_j)(n_j-n_{j-1})+\sum_{j=1}^kd_j(n_{j+1}-n_{j-1}).$$

To give a map of multidegree $\overline d$ is equivalent to specifying a
sequence of quotient bundles
$$ V_{\Bbb P^1}\twoheadrightarrow M_k\twoheadrightarrow \dots \twoheadrightarrow M_1,$$ with 
rank$M_j=n_j$, deg$(M_j)=d_j$, or, by dualizing, a sequence of subbundles
$$S_1\subset\dots\subset S_k\subset V^*_{\Bbb P^1},$$ 
with rank$S_j=n_j$, deg$(S_j)=-d_j$. Let $T_j:=V^*_{\Bbb P^1}/S_{k-j+1}$. The 
Hilbert polynomial of $T_j$ is $P_j(m)=(m+1)(n-n_{k-j+1})+d_{k-j+1}$. 

Let $\Cal H\Cal Q_{\overline d}:=\Cal H\Cal Q_{P_1,\dots ,P_k}
(\Bbb P^1,V^*_{\Bbb P^1})$ be the hyperquot scheme parametrizing flagged
sequences of quotient {\it sheaves} of $V^*_{\Bbb P^1}$, with Hilbert polynomials given by $P_1,\dots ,P_k$.

\proclaim{Theorem 4.1 (\cite{Lau}, \cite{C-F1, C-F2}, \cite{Kim3})} 
$(i)$ $\Cal H\Cal Q_{\overline d}$ is a smooth, irreducible, 
projective variety, of dimension $$\sum_{j=1}^k(n-n_j)(n_j-n_{j-1})+\sum_{j=1}^kd_j(n_{j+1}-n_{j-1}),$$ 
containing $H_{\overline d}$ as an open dense subscheme.

$(ii)$ $\Cal H\Cal Q_{\overline d}$ is a fine moduli space, i.e., there exists
an universal sequence
$$ V^*_{\Bbb P^1\times \Cal H\Cal Q_{\overline d}} 
\twoheadrightarrow \Cal T_k^{\overline d}\twoheadrightarrow \dots 
\twoheadrightarrow 
\Cal T_2^{\overline d}\twoheadrightarrow \Cal T_1^{\overline d}
\twoheadrightarrow 0\tag \dag $$ 
on $\Bbb P^1\times \Cal H\Cal Q_{\overline d}$, such that
each $\Cal T_j^{\overline d}$ is flat over $\Cal H\Cal Q_{\overline d}$, with
relative Hilbert polynomial $P_j(m)$, and having the following property:

For every scheme $X$ over $\Bbb C$, together with a sequence of quotients
$$V^*_{\Bbb P^1\times X}\twoheadrightarrow Q_k\twoheadrightarrow \dots 
\twoheadrightarrow Q_1\tag {\dag\dag}$$
such that each $Q_j$ is flat over $X$, with relative Hilbert polynomial $P_j$,
there exists an unique morphism $\Phi_X:X\rightarrow\Cal H\Cal Q_{\overline d}$
such that the sequence $(\dag\dag)$ is the pull-back of $(\dag)$ via
$(id,\Phi_X)$. 

$(iii)$ Let $\Cal S_j^{\overline d}:=\operatorname{ker}
(V^*_{\Bbb P^1\times \Cal H\Cal Q_{\overline d}}\rightarrow
\Cal T_{k-j+1}^{\overline d})$. Then $\Cal S_j^{\overline d}$ is a
$\underline{vector\ bundle}$ of rank $n_j$ and relative degree $-d_j$ on
$\Bbb P^1\times \Cal H\Cal Q_{\overline d}$. 
\quad\qed\endproclaim

\phantom{X}

\subhead{ 4.2 Generalized Schubert varieties on $H_{\overline d}$ and
$\Cal H\Cal Q_{\overline d}$} \endsubhead

\phantom{X}

The moduli space of maps comes with an universal ``evaluation'' morphism
$$ev: \Bbb P^1\times H_{\overline d} \rightarrow F,$$ given by 
$ ev(t,[\mu])=\mu(t)$, which can be used to pull-back Schubert varieties to 
$H_{\overline d}$. More precisely, 
for $t\in \Bbb P^1$ , $w\in S$, define a subscheme of 
$H_{\overline d}$ by $$\Omega_w(t)=
ev^{-1}(\Omega _w)\bigcap\left ( \{t\} \times  H_{\overline d}\right ).$$ 
(Set-theoretically, $\Omega_w(t)=\{ [\mu]\in 
H_{\overline d}\ |\ \mu(t)\in \Omega_w\}$.)

Alternately, the pull-back $\Omega_w(t)$ of a Schubert variety can be described 
as the degeneracy locus 
$$\{ \operatorname{rank}
(V_p\otimes \Cal O\rightarrow ev^*Q_q)\leq r_w(q,p), 1\leq p\leq n,\ 
q\in\{ n_1,\dots ,n_k\}\}
\bigcap\left ( \{t\} \times  H_{\overline d}\right ),$$
where the $Q_j$'s are the tautological quotient bundles on $F$, and 
$V_1\subset \dots \subset V_{n-1}\subset V_n=V$ 
is our fixed reference flag. This last description may be used to extend
$\Omega_w(t)$ to $\Cal H\Cal Q_{\overline d}$.
 
\definition{ Definition 4.3} $\overline{\Omega}_w(t)$ is the subscheme of 
$\Cal H\Cal Q_{\overline d}$ defined as the degeneracy locus
$$\{ \operatorname{rank}(V_p\otimes \Cal O\rightarrow (\Cal S_q^{\overline d})^*
\leq r_w(q,p), 1\leq p\leq n,\ q\in\{ n_1,\dots ,n_k\}\}
\bigcap\left ( \{t\} \times\Cal H\Cal Q_{\overline d}\right ).\ \square$$
\enddefinition

\phantom{X}

\subhead{4.3 GW-invariants}\endsubhead

\phantom{X}

To define the GW-invariants, we need the following "moving lemma"-type result

\proclaim{ Theorem 4.4 (Moving Lemma)} $(i)$ Let $Y$ be a fixed subvariety of 
$H_{\overline d}$. For any $w\in S$, a corresponding general translate
of $\Omega_w\subset F$, and $t\in \Bbb P^1$, the intersection 
$Y\bigcap\Omega_w(t)$ is either empty, or has pure codimension $\ell (w)$ in $Y$. In particular, for any $w_1,\dots ,w_N\in S$;$\ t_1,\dots ,t_N
\in \Bbb P^1$, and general translates of $\Omega _{w_i}\subset F$, the 
intersection $\bigcap_{i=1}^N \Omega_{w_i}(t_i)$ is either empty, or has pure codimension $\sum_{i=1}^N \ell (w_i)$ in $H_{\overline d}$.

$(ii)$ Moreover, if $t_1,\dots ,t_N$ are $\underline {distinct}$,
then for general 
translates of the $\Omega_{w_i}$ the intersection $\bigcap_{i=1}^N 
\overline {\Omega }_{w_i}(t_i)$ is either empty, or
has pure codimension $\sum_{i=1}^N \ell (w_i)$ 
in $\Cal H\Cal Q_{\overline d}$ and is the Zariski 
closure of $\bigcap_{i=1}^N \Omega_{w_i}(t_i)$.
\endproclaim

\demo{ Proof} $(i)$ This follows from a Theorem of Kleiman (\cite {Kl}), since
$F$ is a homogeneous space.

$(ii)$ The proof will be given in Section 6. \quad\qed
\enddemo

In particular, $\overline\Omega_w(t)$ is the closure of $\Omega_w(t)$ in
$\{ t\}\times\Cal H\Cal Q_{\overline d}$, and has pure codimension $\ell (w)$,
hence, via the identification 
$\{ t\}\times\Cal H\Cal Q_{\overline d}\cong\Cal H\Cal Q_{\overline d}$, 
it determines a class
$[\overline{\Omega}_w(t)]\in H^{2\ell (w)}(\Cal H\Cal Q_{\overline d})$. 
From the above Theorem we get immediately (see e.g. \cite{Be})

\proclaim{Corollary 4.5 } The class $[\overline {\Omega}_w(t)]$ in the cohomology (or Chow) ring of $\Cal H\Cal Q_{\overline d}$ is 
independent of $t\in \Bbb P^1$ and the flag $V_{\bullet }\subset V$. 
\quad\qed
\endproclaim

\proclaim{Corollary 4.6 } Assume that $\sum_{i=1}^N\ell (w_i)=
\operatorname{dim}(H_{\overline d})$. Then, as long as  
$t_1,\dots ,t_N$ are distinct, and we pick general translates of  
the subvarieties $\Omega_{w_i}\subset F$, the number of points in 
$\bigcap_{i=1}^N \Omega_{w_i}(t_i)=
\bigcap_{i=1}^N \overline {\Omega }_{w_i}(t_i)$ 
can be computed as the intersection number
$$\int_{\Cal H\Cal Q_{\overline d}} [\overline\Omega_{w_1}(t_1)]\cup\dots\cup
[\overline\Omega_{w_N}(t_N)]$$ (hence it is 
independent of $t_i$ and the general translates of $\Omega_{w_i}$).\quad\qed
\endproclaim

The corollaries imply that, for general translates of 
$\Omega_{w_1},\dots ,\Omega_{w_N}$ and distinct $t_1,\dots ,t_N$,
we have a well defined intersection number 
$$\langle \Omega_{w_1},\dots ,\Omega_{w_N}\rangle _{\overline d}:
= \cases\operatorname{card}\left (\bigcap_{i=1}^N\Omega_{w_i}(t_i)\right),&
\text{if}\ \sum_{i=1}^N\ell (w_i)=
\operatorname{dim}(H_{\overline d})\\  
0,&\text { otherwise} . \endcases$$
We will call this the {\it Gromov-Witten} invariant associated to the Schubert classes $[\Omega_{w_1}],\dots ,[\Omega_{w_N}]$.

\proclaim{Corollary 4.7} The invariant $I_{3,\overline d}
([\Omega_{w_1}][\Omega_{w_2}][\Omega_{w_3}])$ defined using Kontsevich's space of stable maps $\overline M_{0,3}(F,\overline d )$ coincides with
$\langle \Omega_{w_1},\Omega_{w_2},\Omega_{w_3}\rangle _{\overline d}$.
\quad\qed\endproclaim

\remark{\bf Remark 4.8}  By the preceding Corollary we have
$$[\Omega_u]*[\Omega_v]=\sum_{\overline d}q_1^{d_1}\dots q_k^{d_k}
\sum_{w\in S}\langle\Omega_u,\Omega_v,\Omega_w\rangle _{\overline d}
[\Omega_{\check w}].$$ 
Note that we have also
$$[\Omega_{w_1}]*[\Omega_{w_2}]*\dots *[\Omega_{w_N}]
=\sum_{\overline d}q_1^{d_1}\dots q_k^{d_k}
\sum_{w\in S}\langle \Omega_{w_1},\dots ,\Omega_{w_N},\Omega_w
\rangle _{\overline d}[\Omega_{\check w}].\tag 4.1$$
We will use (4.1) in Section 6, for the proof of Theorem 3.16. \quad\qed
\endremark

\phantom{X}

Finally, we record for later use a generalization of (part of) Theorem 4.4, due to B. Kim (\cite{Kim3, Corollary 3.2}). 

For every irreducible closed subvariety $Y\subset F$ and every
$t\in \Bbb P^1$, we denote by $Y(t)$ the preimage 
$ev^{-1}(Y)\bigcap\{ y\}\times H_{\overline d}$ and by $\overline{Y(t)}$ 
the closure of $Y(t)$ in $\{ y\}\times \Cal H\Cal Q_{\overline d}$.

\proclaim{Proposition 4.9} Let $Y_1,\dots ,Y_N$ be closed, irreducible 
subvarieties in $F$, and let $t_1,\dots ,t_N$ be distinct points in $\Bbb P^1$.
Assume that $\sum_{i=1}^N\operatorname{codim}Y_i=\operatorname{dim}H_{\overline d}$. Then for general translates of $Y_i$, the intersection scheme
$\bigcap_{i=1}^NY_i(t_i)$ is either empty, or it consists of finitely many 
reduced points. Moreover,
$$\bigcap_{i=1}^NY_i(t_i)=\bigcap_{i=1}^N\overline{Y_i(t_i)},$$ and the
cardinality of this set is equal to the intersection number
$$\int_{\Cal H\Cal Q_{\overline d}}[\overline{Y_1(t_1)}]\cup\dots\cup
[\overline{Y_N(t_N)}].\quad\qed $$
\endproclaim

\heading{5. The boundary of $\Cal H\Cal Q_{\overline d}$}\endheading

The space $H_{\overline d}$ is the largest subscheme  of 
$\Cal H\Cal Q_{\overline d}$ with the property 
that on $\Bbb P^1\times H_{\overline d}$ the sheaf injections in the universal sequence
$$ 0\hookrightarrow 
\Cal S_1^{\overline d}\hookrightarrow \Cal S_2^{\overline d}\hookrightarrow 
\dots \hookrightarrow 
\Cal S_k^{\overline d}\hookrightarrow V^*_{\Bbb P^1\times \Cal H\Cal Q_{\overline d}}$$
are vector bundle inclusions. The boundary of $\Cal H\Cal Q_{\overline d}$,
by which we mean the complement of $H_{\overline d}$, is therefore the locus
$\Cal B_{\overline d}$ such that on $\Bbb P^1\times\Cal B_{\overline d}$ 
some of these maps degenerate. 
In this section we will study the restrictions of the generalized Schubert varieties $\overline {\Omega}_w(t)$
to $\Cal B_{\overline d}$. We start with a description of $\Cal B_{\overline d}$
itself, for which the following construction (taken from \cite{C-F2, 2.2}) is
needed:

Let $\overline e=(e_1,\dots ,e_k)$ be a multiindex satisfying the
conditions 

$(5.1)\ \ e_i\leq \operatorname {min}(n_i,d_i),\ \text {for}\ 1\leq i\leq k$,

$(5.2)\ \ e_i-e_{i-1}\leq n_i-n_{i-1},\ \text {for}\ 2\leq i\leq k$.

\flushpar(cf. \cite {C-F2, Lemma 2.1}).

For each $1\leq i\leq k$, let $\pi_i:\Cal G_i\longrightarrow
\Bbb P^1\times \Cal H\Cal Q_{\overline d-\overline e}$ be the Grassmann bundle of $e_i$-dimensional quotients of $\Cal  S_i^{\overline d-\overline e}$. Let
$\Cal G_{\overline e}$ 
be the fiber product of the $\Cal G_i$'s over $\Bbb P^1\times 
\Cal H\Cal Q_{\overline d-\overline e}$, with projection
$\pi :\Cal G_{\overline e}
\longrightarrow \Bbb P^1\times \Cal H\Cal Q_{\overline d-\overline e}$.

For each $1\leq i\leq k$, let
$$0\longrightarrow K_i\longrightarrow\pi _i^*\Cal S_i^{\overline d-\overline e}
\longrightarrow L_i\longrightarrow 0$$
be the universal sequence on $\Cal G_i$. $K_i$ and $L_i$ are vector bundles, of ranks
$n_i-e_i$ and $e_i$ respectively. On $\Cal G_{\overline e}$ we have 
the following diagram:
$$\minCDarrowwidth{7 mm}
\CD
0@. ... @. 0 @. 0 @. ... @. 0 @. \\
@VVV @. @VVV @VVV @. @VVV @.\\
K_1@. ... @. K_i @. K_{i+1}
@. ... @. K_k @. \\
@VVV @. @VVV @VVV @. @VVV @. \\ 
\pi ^*\Cal S_1^{\overline d-\overline e} @>>>
  ... @>>> \pi ^*\Cal S_i^{\overline d-\overline e} 
@>>> \pi ^*\Cal S_{i+1}^{\overline d-\overline e} @>>> ... @>>> 
\pi ^*\Cal S_k^{\overline d-\overline e} @>>> V^*_{\Cal G_{\overline e}}\\
@VVV @. @VVV @VVV @. @VVV @. \\
L_1 @.  ... @. L_i @. L_{i+1} @. ...  @. 
L_k @. \\
@VVV @. @VVV @VVV @. @VVV @. \\
0@. ... @. 0 @. 0 @. ... @. 0 @. \endCD$$

Let ${\Cal U}_{\overline e}$ be the locally-closed subscheme of
$\Cal G_{\overline e}$ determined by the closed conditions
$$\operatorname {rank}(K_i\longrightarrow L_{i+1})=0,\ \ \ 
\text {for}\ i=1,\dots ,k-1,\tag 5.3$$
and the open conditions
$$\operatorname {rank}(K_i\longrightarrow 
V^*_{\Cal G_{\overline e}})=n_i-e_i,\ \ \ \text {for}\ i=1,\dots ,k.\tag 5.4 $$

\proclaim{Lemma 5.1} ${\Cal U}_{\overline e}$ is smooth, irreducible, of 
dimension 
$$1+\operatorname{dim}(\Cal H\Cal Q_{\overline d})-\sum_{i=1}^ke_i(n_{i+1}-n_i)
-\sum_{i=1}^ke_i(e_i-e_{i-1}).$$
The projection $\pi :\Cal U_{\overline e}\rightarrow 
\Bbb P^1\times\Cal H\Cal Q_{\overline d-\overline e}$ 
is smooth, and its image contains $\Bbb P^1\times H_{\overline d-\overline e}$.
\endproclaim

\demo{Proof} For a vector bundle $E$ on a scheme $X$, we denote by $G_e(E)$ the Grassmann bundle of $e$-dimensional quotients of $E$, for some $0\leq e\leq\operatorname{rank}(E)$.

Let $\Cal V_{\overline e}\subset\Cal G_{\overline e}$ be
the open subscheme defined by the conditions (5.4) and put 
$\Cal V:=\pi (\Cal V_{\overline e})$. Obviously, $\Cal V$ is open in 
$\Bbb P^1\times\Cal H\Cal Q_{\overline d-\overline e}$, and contains
$\Bbb P^1\times H_{\overline d-\overline e}$. 
The Lemma is a consequence of the observation that 
$\Cal U_{\overline e}$ can be constructed as a sequence of $k$ Grassmann bundles
over $\Cal V$ as follows:

$\bullet$ start with 
$\rho_1:G_{e_1}(\Cal S_1^{\overline d-\overline e})\rightarrow\Cal V$ 
with universal subbundle $K_1$;

$\bullet$ next, form $\rho_2:G_{e_2}(\rho_1^*\Cal S_2^{\overline d-\overline e}/K_1)\rightarrow G_{e_1}(\Cal S_1^{\overline d-\overline e})$
with universal subbundle $\Cal L_2$, and let $K_2$ be the natural induced
extension
$$0\rightarrow \rho_2^*K_1\rightarrow K_2\rightarrow\Cal L_2\rightarrow 0;$$

$\bullet$ continue by forming 
$\rho_3:G_{e_3}(\rho_2^*\rho_1^*\Cal S_3^{\overline d-\overline e}/K_2)
\rightarrow G_{e_2}(\rho_1^*\Cal S_2^{\overline d-\overline e}/K_1)$, 
with universal subbundle $\Cal L_3$, and let $K_3$ be the natural extension
$$0\rightarrow \rho_3^*K_2\rightarrow K_3\rightarrow\Cal L_3\rightarrow 0,$$
and so on.

We will use this description of $\Cal U_{\overline e}$ in the proof of
Quantum Pieri (see 6.3).
\quad\qed
\enddemo

\proclaim{Theorem 5.2}There exist morphisms $\phi_{\overline e}\ :\ 
\Cal U_{\overline e}\longrightarrow 
\Cal H\Cal Q_{\overline d}$ such that

$(i)$ If $\operatorname {rank}_{(t,x)}\Cal T_{k-i+1}^{\overline d}= 
n-n_i+e_i$ 
for every $1\leq i\leq k$ at a point 
$(t,x)\in\Bbb P^1\times\Cal H\Cal Q_{\overline d}$, then 
$x\in \phi_{\overline e}(\Cal U_{\overline e})$. 
In particular $\Cal B_{\overline d}$ is covered by
the union of $\phi_{\overline e}(\Cal U_{\overline e})$, where $\overline e$
ranges over all (nonzero) multiindices satisfying (5.1) and (5.2).

$(ii)$ The restriction of $\phi_{\overline e}$ to
$\pi ^{-1}(\Bbb P^1\times H_{\overline d-\overline e})$ is an isomorphism onto its image.
\endproclaim

\demo{Proof} See \cite{C-F2, Theorem 2.3 $(ii)$}.\quad\qed
\enddemo

\proclaim {Lemma 5.3} 
$$\phi_{\overline e}^{-1}(\overline{\Omega}_w(t))
=\pi ^{-1}(\Bbb P^1\times\overline{\Omega}_w(t))
\bigcup \tilde\Omega _w^{\overline e}(t),$$
$\tilde\Omega _w^{\overline e}(t)$ being the locus inside 
$\Cal U_{\overline e}(t):=\pi ^{-1}(\{ t\}\times 
\Cal H\Cal Q_{\overline d-\overline e})$ where 
$$\operatorname {rank}\ (V_p\otimes \Cal O\rightarrow K_q^*)
\leq r_w(q,p),\tag 5.5$$ for all $p=1,\dots ,n$, $q\in\{ n_1,\dots ,n_k\} $.
\endproclaim

\demo{Proof} See \cite{C-F2, Lemma 3.1}.\quad\qed
\enddemo

Following \cite{C-F2, Section 3} we will describe now the locus
$\tilde\Omega _w^{\overline e}(t)$ of Lemma 5.3. The analysis there can be 
applied in our case without any changes and the only reason for reproducing 
part of it here is to fix the somewhat elaborate notation needed.

Let $a:=\operatorname {card}\{n_i-e_i\ \mid \ i=1,\dots ,k\}$. 
Set $e_0=e_{k+1}=0$ and define a partition of $[0,k+1]$ as follows:
$$i_0=0,\ i_j=\operatorname {min}\{\ i\ \mid 
\ n_i-e_i\geq n_{i_{j-1}}-e_{i_{j-1}}+1\},
\ \operatorname {for}\ 1\leq j\leq a,\ i_{a+1}=k+1.$$
Let $m_j=n_{i_j}-e_{i_j}$, for $j=0,1,\dots ,a$.
By definition, on each of the intervals $$[1,i_1-1],\
[i_1,i_2-1],\dots ,[i_{a-1},i_a-1],\ [i_a,k]$$ $n_i-e_i$ is constant, equal 
respectively to $0,m_1,\dots ,m_a$, and the corresponding bundles 
$K_i^*$ are all isomorphic.
Therefore we can restrict the set of rank conditions (5.5), defining 
$\tilde\Omega _w^{\overline e}(t)$ inside $\Cal U_{\overline e}(t)$, to
$$\operatorname {rank}\ (V_p\otimes \Cal O\rightarrow K_q^*)
\leq r_w(q,p),\ 1\leq p\leq n,\ q\in\{ n_{i_1},\dots ,n_{i_a}\} .\tag 5.6$$

Moreover, we can further modify (5.6).  
Define recursively $\bold{r}:=(r_{j,p})_{1\leq p\leq n,\ 1\leq j\leq a}$ as follows:
$$r_{1,p}=\operatorname{min}\{ r_w(n_{i_1},p),\ m_1\} ,\ 1\leq p\leq n,$$
$$r_{j,p}=\operatorname{min}\{ r_w(n_{i_j},p),\ r_{j-1,p}+m_j-m_{j-1}\} ,
\ 1\leq p\leq n,\ 2\leq j\leq a.$$

\proclaim{Lemma 5.4} The conditions
$$\operatorname {rank}\ (V_p\otimes \Cal O\rightarrow K_{i_j}^*)
\leq r_{j,p},\ 1\leq p\leq n,\ 1\leq j\leq a\tag 5.7$$ 
define the same degeneracy
locus $\tilde\Omega _w^{\overline e}(t)$ in $\Cal U_{\overline e}(t)$.
\endproclaim

\demo{ Proof } See \cite{C-F2, 3.5}.\quad\qed\enddemo

\proclaim{Lemma 5.5} $(i)$ There exists a (unique) permutation
$\tilde w^{\overline e}\in S_n$ such that if
$\tilde w^{\overline e}(q)>\tilde w^{\overline e}(q+1)$, then 
$q\in\{ m_1,\dots ,m_a\} $, and
$r_{j,p}=r_{\tilde w^{\overline e}}(m_j,p)$, for all
$1\leq p\leq n,\ 1\leq j\leq a$.

$(ii)$ $\ell (w)-\ell (\tilde w^{\overline e})\leq\sum_{i=1}^ke_i(n_{i+1}-n_i)$.
\endproclaim

\demo{ Proof } The following explicit construction of
$\tilde w^{\overline e}$ is taken from \cite{C-F2, Lemma 3.6}.

For each $j=0,1,\dots ,a+1$, define sets $W_j(w)$, by
$$W_0(w)=\emptyset, \ W_j(w)=\{w(1),\dots ,w(n_{i_j})\}.$$
Also we define sets $Z_j(w)$, and {\it ordered} sets $\tilde Z_j(w)$, such that

(a) $\operatorname {card}Z_j(w)=n_{i_j}-m_{j-1}$, 

(b) $\operatorname {card}\tilde Z_j(w)=m_j-m_{j-1}$,

(c) $\tilde Z_j(w)\bigcap\tilde Z_{j'}(w)=\emptyset$ if $j\neq j'$,
 
(d) $\bigcup_{j=1}^{a+1}\tilde Z_j(w)=[1,n],$
  
\flushpar by the following recursive procedure:

Let $Z_0(w)=\tilde Z_0(w)=\emptyset$. If $\tilde Z_i(w)$ has been already
defined for $i=0,1,\dots ,j-1$, let
$$Z_j(w)=W_j(w)\setminus \left(\bigcup _{i=1}^{j-1}\tilde Z_i(w)\right ).$$
Arrange $Z_j(w)$ in increasing order
$$Z_j(w)=\{z_{j,1}<\dots <z_{j,n_{i_j}-m_{j-1}}\},$$
and set 
$$\tilde z_{m_{j-1}+1}:=z_{j,1},\ \tilde z_{m_{j-1}+2}:=z_{j,2},\ \dots ,
\ \tilde z_{m_j}:=z_{j,m_j-m_{j-1}},$$
$$\tilde Z_j(w):=\{\tilde z_{m_{j-1}+1}<\tilde z_{m_{j-1}+2}< \dots <
\tilde z_{m_j}\}.$$
Now define $\tilde w^{\overline e}(q)=\tilde z_q$,
for all $1\leq q\leq n$.  

The estimate $(ii)$ follows (cf. \cite{C-F2, Lemma 3.8})
by noticing first that
the difference $\ell (w)-\ell (\tilde w^{\overline e})$ is maximized
by the longest permutation $w^{\circ}$, defined by $w^{\circ}(i)=n-n_j+i-n_{j-1}$, for all $n_{j-1}+1\leq i\leq n_j$, 
$1\leq j\leq k+1$. For this case one computes directly that we have in 
fact the equality
$$\ell (w^{\circ})-\ell ((\tilde {w^{\circ}})^{\overline e})=\sum_{i=1}^ke_i(n_{i+1}-n_i).\quad\qed$$\enddemo

\remark{\bf Remark 5.6} $(i)$ In the terminology of \cite{F1}, 
Lemma 5.5 $(i)$ says that $\bold{r}$ is a {\it permissible} collection of rank
numbers.

$(ii)$ Let $\tilde F_{\overline e}:=F(m_1,\dots ,m_a, V)$ be
the partial flag variety corresponding to the $m_i$'s.
The sequence of quotient bundles
$$V\otimes \Cal O_{\Cal U_{\overline e}(t)}\twoheadrightarrow K^*_{n_{i_a}}
\twoheadrightarrow \dots\twoheadrightarrow K^*_{n_{i_1}}$$ is the pull-back via 
an uniquely determined morphism $\psi_{\overline e}(t):\Cal U_{\overline e}(t)
\rightarrow\tilde F_{\overline e}$ of the tautological sequence on
$\tilde F_{\overline e}$. By Lemmas 5.4 and 5.5 $(i)$,
$\tilde w^{\overline e}$ defines a
Schubert variety on $\tilde F_{\overline e}$, and we have
$\tilde\Omega _w^{\overline e}(t)=\psi_{\overline e}(t)^{-1}
(\Omega_{\tilde w^{\overline e}})$. 
\quad\qed\endremark

Finally, we spell out in more detail what the analysis in this section
says for some special cases.

\proclaim{ Lemma 5.7 } Let $(e_1,\dots ,e_k)$ be a multiindex. 
Fix $1\leq j\leq k$ and $1\leq i\leq n_j$ and consider the
cycle $\alpha_{i,j}:=s_{n_j-i+1}\cdot\dots\cdot s_{n_j}$. Then
$$\tilde\alpha_{i,j}^{\overline e}=
\cases{\alpha_{i,j},\ \ \text{if}\ \ e_j=0,}\\
{\alpha_{i,j}\cdot s_{n_j}\cdot\dots\cdot s_{n_j-e_j+1},\ \ \text{if}\ \
1\leq e_j<i,}\\{id,\ \ \text{if}\ \ i\leq e_j.}
\endcases  $$ In particular, $\ell (\alpha_{i,j})-\ell (\tilde\alpha_{i,j}^{\overline e})\leq e_j$, with equality iff $e_j\leq i$.
\endproclaim

\demo{ Proof} Immediate from the construction of 
$\tilde\alpha_{i,j}^{\overline e}$ in Lemma 5.5. \quad\qed\enddemo

\proclaim{ Lemma 5.8 (cf. \cite{C-F2, Lemma 3.9})}  Assume in addition that 
$\overline e\neq(0\dots ,0)$. Then

$(i)$ $ \sum_{i=1}^ke_i(e_i-e_{i-1})\geq e_j$, for $1\leq j\leq k$. In
particular $\sum_{i=1}^ke_i(e_i-e_{i-1})\geq 1$.

$(ii)$ There exists $1\leq j\leq k$ such that 
$\sum_{i=1}^ke_i(e_i-e_{i-1})=e_j$ iff the following holds:

there are integers $1\leq h_1<h_2<\dots <h_m\leq j\leq l_m<\dots <l_2<l_1\leq k$
such that $$e_i=\cases{0,\ \text{for}\ i\in [1,h_1-1]\cup [l_1+1,k],}\\
{1,\ \text{for}\ i\in [h_1,h_2-1]\cup [l_2+1,l_1],}\\
{2,\ \text{for}\ i\in [h_2,h_3-1]\cup [l_3+1,l_2],}\\
{\dots}\\
{m,\ \text{for}\ i\in [h_m,l_m]}\endcases$$ 
(in particular, $e_j=m$).

$(iii)$ Let $\overline e_{\bold{h}\bold{l}}$ denote a multiindex as in 
$(ii)$, and let
$w\in S$ be any permutation. Let $\tilde w^{\overline e_{\bold{h}\bold{l}}}$ be the permutation given by Lemma 5.5 $(i)$. Then
$$\ell (w)-\ell(\tilde w^{\overline e_{\bold{h}\bold{l}}})=\sum_{i=1}^ke_i(n_{i+1}-n_i)=
\sum_{c=1}^m(n_{l_c+1}-n_{h_c})$$
iff for {\it every} $1\leq i\leq m$ we have
$$w(n_{h_i})>\operatorname{max}\{w(n_{h_i}+1),
\dots ,w(n_{l_i}),w(n_{l_i+1})\}.$$ 
In this case, $\tilde w^{\overline e_{\bold{h}\bold{l}}}=w\cdot\gamma_{h_m,l_m}\cdot
\gamma_{h_{m-1},l_{m-1}}\cdot\dots\cdot\gamma_{h_1,l_1}$, where
$\gamma_{h,l}$ denotes the cyclic permutation $s_{n_h}\cdot\dots\cdot
s_{n_{l+1}-1}$ (cf. the paragraph before Theorem 3.1).

\endproclaim

\demo{ Proof}  $(i)$ First, using the easy identity
$$\sum_{i=1}^ke_i(e_i-e_{i-1})={1\over 2}[e_1^2+(e_2-e_1)^2+\dots +
(e_k-e_{k-1})^2+e_k^2],$$
and the change of variables 
$$x_1=e_1,\ x_2=e_2-e_1,\ \dots ,\ x_k=e_k-e_{k-1},\ x_{k+1}=e_k,$$
the inequality in $(i)$ becomes
$$\sum_{i=1}^j(x_i^2-2x_i)+\sum_{i=j+1}^{k+1}x_i^2\geq 0,\tag 5.8$$
with the additional constraint $x_{k+1}=\sum_{i=1}^kx_i$. Now (5.8) is equivalent to
$$\sum_{i=1}^j(x_i-1)^2+\sum_{i=j+1}^{k+1}x_i^2\geq j.\tag 5.9$$
Making another change of variables 
$$y_1=x_1-1,\dots ,y_j=x_j-1,y_{j+1}=x_{j+1},\dots ,y_{k+1}=x_{k+1},$$ 
we are reduced to proving
$$\sum_{i=1}^{k+1}y_i^2\geq j,\tag 5.10$$
subject to the constraint
$$y_{k+1}=j+\sum_{i=1}^ky_i.\tag 5.11$$
Replacing in (5.10) $j$ by $y_{k+1}-\sum_{i=1}^ky_i$, we get
$$(y_{k+1}^2-y_{k+1})+\sum_{i=1}^k(y_i^2+y_i)\geq 0.\tag 5.12$$
Since $y_i$, $1\leq i\leq k+1$ are integers, each paranthesis in (5.12) is
nonnegative. This proves $(i)$.

$(ii)$ We have equality in (5.12) iff $y_{k+1}$ is equal to either $0$ or $1$,
and each $y_i,\ 1\leq i\leq k$ is equal to either $0$ or $-1$. Using (5.11), we
see that equality occurs in one of the following two cases

$\bullet$ either $y_{k+1}=0$, exactly $j$ among $y_1\dots ,y_k$ are equal to $-1$, and the rest of them are equal to $0$,

$\bullet$ or $y_{k+1}=1$, exactly $j-1$ among $y_1\dots ,y_k$ are equal to $-1$, and the rest of them are equal to $0$.

Changing the variables back to $e_i$, the statement in $(ii)$ is obtained.

$(iii)$ follows by the construction of $\tilde w^{\overline e_{\bold{h}\bold{l}}}$. 
\quad\qed\enddemo

\heading{6. Proofs of the Moving Lemma, the Quantum Pieri Formula, and the Orthogonality Theorem}\endheading

This section is devoted to the proofs of Theorem 4.4 $(ii)$, Theorem 3.1, and
Theorem 3.16. For this purpose we will use heavily the structure of the 
boundary of $\Cal H\Cal Q_{\overline d}$, described in the preceding section.

Throughout the rest of the paper, we will work with suitable 
general translates of the Schubert varieties on $F$.

\phantom{X}

\subhead{6.1 Proof of Theorem 4.4 $(ii)$ 
(cf. \cite{C-F2, Theorem 4.1})}\endsubhead

\phantom{X}

We proceed by induction on $\overline d$.
If $d_1=\dots =d_k=0$, then $\Cal H\Cal Q_{\overline d}=H_{\overline d} =F$ and there is nothing to prove. 
Assume that the statement is true for all $\overline f$ such that $f_i\leq d_i$, $1\leq i\leq k$ and $f_j<d_j$ for some $1\leq j\leq k$.
Let $c:=\sum_{i=1}^N \ell (w_i)$. By $(i)$ and Theorem 5.2 $(i)$, 
it is enough to show that
$$\operatorname{codim}_{\Cal H\Cal Q_{\overline d}}\left (\bigcap _{i=1}^N \overline {\Omega}_{w_i}(t_i)\right )\bigcap\ \phi_{\overline e}
(\Cal U_{\overline e})>c,\tag 6.1$$
for every multiindex $\overline e\neq(0,\dots ,0)$, satisfying the conditions
(5.1) and (5.2). Using Theorem 5.2 $(ii)$ and Lemma 5.1, the inequality (6.1)
will follow if we prove that the codimension of
$\bigcap_{i=1}^N \phi_{\overline e}^{-1}
\left (\overline {\Omega}_{w_i}(t_i)\right )$ 
in $\Cal U_{\overline e}$ is greater than
$$c-(\operatorname {dim}\ \Cal H\Cal Q_{\overline d}-\operatorname {dim}\ 
\Cal U_{\overline e})=c+1-\sum_{i=1}^ke_i(n_{i+1}-n_i)-
\sum_{i=1}^ke_i(e_i-e_{i-1}).$$
By Lemma 5.3, we have to prove the same estimate for the codimension of
$$\bigcap_{i=1}^N\left (\pi ^{-1}(\Bbb P^1\times \overline{\Omega }_{w_i}(t_i))
\bigcup \ {\tilde \Omega}^{\overline e}_{w_i}(t_i)\right )\tag 6.2$$
in $\Cal U_{\overline e}$. Since the points $t_1,\dots ,t_N$ are distinct, the only possibly nonempty intersections in (6.2) contain either no
${\tilde \Omega}^{\overline e}_{w_i}(t_i)$ term, or only one such term. If 
there is no such term, the required inequality follows from the induction
assumption on $\Cal H\Cal Q_{\overline d-\overline e}$ and the fact that
$\pi $ is a smooth map. After possibly renumbering the points $t_i$,
to finish the proof it suffices to show the estimate for
$$W\bigcap\ {\tilde \Omega}^{\overline e}_{w_N}(t_N)\subset
\Cal U_{\overline e}(t_N), \tag 6.3$$
where
$$W:=\bigcap_{i=1}^{N-1}\left (\pi^{-1}(\{ t_N\}\times \overline
{\Omega }_{w_i}(t_i))\right ),$$ and
$\Cal U_{\overline e}(t_N)=\pi ^{-1}(\{ t_N\}\times\Cal H\Cal Q_
{\overline d-\overline e})$. 
By the induction assumption,
W has codimension $c-\ell (w_N)$ in $\Cal U_{\overline e}(t_N)$,
while by Remark 5.6 $(ii)$ and Kleiman's Theorem on transversality of general
translates, the intersection (6.3) has codimension 
$\ell (\tilde w^{\overline e})$ in $W$. The estimate follows now from
Lemma 5.5 $(ii)$ and Lemma 5.8 $(i)$. \quad\qed

\phantom{X}

\subhead{6.2 Computing GW-invariants via degenerations}\endsubhead

\phantom{X}

For the proofs of Quantum Pieri and Orthogonality, we need to compute certain
invariants $\langle \Omega_{w_1},\dots ,\Omega_{w_N}\rangle _{\overline d}$.
The technique we will use is to degenerate the intersection
$\bigcap_{i=1}^N \overline\Omega_{w_i}(t_i)$ by allowing some of the points
$t_i$ to coincide. This procedure may lead to contributions supported on
the boundary, which can be evaluated using the analysis in Section 5.
In this subsection we summarize some results of this type.

The following is an easy consequence of Proposition 4.9 and Theorem 4.4. 
For a proof, see for instance \cite{Be, Lemma 2.5}.

\proclaim{Lemma 6.1} Let $Y_1$, $Y_2$ be subvarieties in $F$ such that
$\operatorname{codim}Y_1+\operatorname{codim}Y_1=\operatorname{dim}
\Cal H\Cal Q_{\overline d}$, and let $t_1,t_2\in \Bbb P^1$ be distinct points.
Assume $\overline d\neq(0,\dots ,0)$. Then
$$\int_{\Cal H\Cal Q_{\overline d}}[\overline{Y_1(t_1)}]\cup
[\overline{Y_2(t_2)}]=0.$$
In particular, for any $v,w\in S$,
$$\langle\Omega_v,\Omega_w\rangle_{\overline d}=
\cases{1,\ \ \text{if}\ \overline d=(0,\dots ,0)\ \text{and}\ v=\check w}\\
{0,\ \ \text{otherwise}}\endcases . \quad\qed$$\endproclaim

\proclaim{Lemma 6.2} Let ${\overline d}$ be a multiindex, and let
$v_1,\dots ,v_N,w_1,\dots ,w_M\in S$ satisfy $\sum_{i=1}^N\ell (v_i)+\sum_{j=1}^M\ell (w_j)=\operatorname{dim}H_{\overline d}$. Let
$y,t_1,\dots ,t_M\in \Bbb P^1$ be distinct points. 
Assume that the conclusion of 
Theorem 4.4 $(ii)$ holds for the intersection
$$\overline\Omega_{v_1}(y)\bigcap\dots\bigcap\overline\Omega_{v_N}(y)\bigcap
\overline\Omega_{w_1}(t_1)\bigcap\dots\bigcap\overline\Omega_{w_M}(t_M).
\tag 6.4$$ Let $Y:=\Omega_{v_1}\bigcap\dots\bigcap\Omega_{v_N}\subset F$.
Then $$\langle\Omega_{v_1},\dots ,\Omega_{v_N}, 
\Omega_{w_1},\dots ,\Omega_{w_M}\rangle_{\overline d}=\int_{\Cal H\Cal Q
_{\overline d}}[\overline{Y(y)}]\cup [\overline\Omega_{w_1}(t_1)]\cup\dots 
\cup[\overline\Omega_{w_N}(t_N)] .$$
\endproclaim

\demo{Proof} By Theorem 4.4 $(i)$, the intersection 
$$\Omega_{v_1}(y)\bigcap\dots\bigcap\Omega_{v_N}(y)\bigcap
\Omega_{w_1}(t_1)\bigcap\dots\bigcap\Omega_{w_M}(t_M)\subset
H_{\overline d}$$ has pure dimension $0$, and by assumption it coincides with
(6.4). On the other hand, we have
$$\Omega_{v_1}(y)\bigcap\dots\bigcap\Omega_{v_N}(y)=ev^{-1}(Y)
\bigcap\{ y\}\times H_{\overline d}.$$
The Lemma follows now from the definition of 
$\langle\dots\rangle _{\overline d}$,
Corollaries 4.5 and 4.6, and Proposition 4.9. \quad\qed\enddemo

\proclaim{Proposition 6.3} Let ${\overline d}\neq (0,\dots ,0)$ be a multiindex, let $u,v,w\in S$ be such that $\ell (u)+\ell (v)+\ell (w)=
\operatorname{dim}H_{\overline d}$, and let
$y,t\in \Bbb P^1$ be distinct points. Denote
$$Z:=\overline\Omega_u(y)\bigcap\overline\Omega_v(y)\bigcap
\overline\Omega_w(t)\subset\Cal H\Cal Q_{\overline d}.$$ Assume that $Z$ is
either empty, or purely 0- dimensional. Then

$(i)$ $Z$ is contained in $\Cal B_{\overline d}$.

$(ii)$ $[Z]=[\overline\Omega_u(y)]\cup [\overline\Omega_v(y)]\cup
[\overline\Omega_w(t)]$ is a cycle of degree $\langle\Omega_u,\Omega_v,\Omega_w
\rangle _{\overline d}$. 
\endproclaim

\demo{Proof} $(i)$ Write
$\Omega_u\bigcap\Omega_v=Y$ inside $F$.
Let $Z'$ be the (largest) subscheme of $Z$ supported on
$H_{\overline d}$. Then, as in the proof of Lemma 6.2, we have
$$Z'=Y(y)\bigcap\Omega_w(t),$$
hence, by Proposition 4.9,
$$\operatorname{card}Z'=\int_{\Cal H\Cal Q
_{\overline d}}[\overline{Y(y)}]\cup[\overline\Omega_w(t)].$$ 
By Lemma 6.1, $Z'$ is empty.

$(ii)$ Let $V:=\overline\Omega_v(y)\bigcap
\overline\Omega_w(t)\subset\Cal H\Cal Q_{\overline d}$, and consider the trivial
family $\Bbb P^1\times V\subset\Bbb P^1\times\Cal H\Cal Q_{\overline d}$ over $\Bbb P^1$. Let $\rho :X\hookrightarrow\Bbb P^1\times\Cal H\Cal Q_{\overline d}
\overset{pr}\to{\longrightarrow}\Bbb P^1$ be the family whose fibre over
$x\in \Bbb P^1$ is the generalized Schubert variety $\overline\Omega_u(x)$.
It follows from Theorem 4.4 that $\rho$ is a fibre bundle map (see e.g. \cite{Be, Corollary 2.4}), and in particular it is flat. Since the intersection
$\left (\Bbb P^1\times V\right )\bigcap X$ 
is obviously proper over $\Bbb P^1$, the Proposition follows from Example 10.2.1 in \cite{F2}. \quad\qed\enddemo

\subhead{6.3 Proof of Quantum Pieri}\endsubhead

\phantom{X}

We formulate first an auxiliary Lemma, for which we introduce some notation.

Let $X$ be a scheme. Let $V$ be an $n$-dimensional complex vector
space and let $B_i$, $1\leq i\leq k$ be vector bundles on $X$, of ranks $b_i$ respectively. Fix $1\leq j\leq k$.
Assume that we are given  a sequence of generically injective maps
$$B_1\rightarrow \dots\rightarrow B_{j-1}\rightarrow B_j\rightarrow B_{j+1}
\rightarrow\dots\rightarrow B_k\rightarrow V^*\otimes\Cal O_X.$$
Moreover, assume that $B_i\rightarrow V^*$ is an injective map of bundles
for $1\leq i\leq j$.
Fix $0\leq e\leq b_j-b_{j-1}$, 
and let $\rho :G_e(B_j/B_{j-1})\rightarrow X$ be the Grassmann bundle of $e$-dimensional quotients of $B_j/B_{j-1}$, 
with universal sequence
$$0\rightarrow L\rightarrow \rho^*(B_j/B_{j-1})\rightarrow Q\rightarrow 0.$$
Let $K$ be the natural induced extension
$$0\rightarrow \rho^*B_{j-1}\rightarrow K\rightarrow L\rightarrow 0,$$
i.e., $K$ is the kernel of $\rho^*B_j\rightarrow Q$. 

Let $V_1\subset\dots
\subset V_{n-1}\subset V$ be a fixed flag, and let $w\in S_n$ be a permutation
such that if $w(q)>w(q+1)$, then $q\in\{ b_1, \dots ,b_{j-1}, b_j-e,b_{j+1},
\dots ,b_k\} $.
Denote by $\bold{D}_w$ the degeneracy locus on $G_e(B_j/B_{j-1})$ determined by
$$\operatorname{rank}(V_p\otimes\Cal O\rightarrow (\rho^*B_q)^*)\leq r_w(q,p),
\ 1\leq p\leq n,\ q\in\{ b_1,\dots ,b_{j-1},b_{j+1},\dots ,b_k\} ,\tag 6.5$$ 
and
$$ \operatorname{rank}(V_p\otimes\Cal O\rightarrow K^*)\leq r_w(b_j-e,p), 
\ 1\leq p\leq n.\tag 6.6$$
Define a permutation $\hat w\in S_n$ as follows:

$\bullet$ let $\{ z_1<z_2<\dots <z_{b_j-b_{j-1}}\} $ be the set
$\{ w(b_{j-1}+1), w(b_{j-1}+2),\dots ,w(b_j)\} $, ordered increasingly;

$\bullet$ if $q\notin\{ b_{j-1}+1, b_{j-1}+2,\dots ,b_j\} $, set
$\hat w(q)=w(q)$;

$\bullet$ for $1\leq i\leq b_j-b_{j-1}$, set $\hat w(b_{j-1}+i)=z_i$.

\proclaim{Lemma 6.4} $(i)$ 
The image $\rho (\bold{D}_w)$ is the degeneracy locus $\bold{D}_{\hat w}$ on $X$ defined by
$$\operatorname{rank}(V_p\otimes\Cal O\rightarrow B_q^*)\leq r_{\hat w}(q,p), 
\ 1\leq p\leq n,\ q\in\{ b_1, \dots  ,b_k\} .\tag 6.7$$

$(ii)$ The restriction of $\rho $ to $\bold{D}_w$ has positive dimensional
fibres, unless 
$$w(b_{j-1}+1)>w(b_j),\tag 6.8$$ in which case 
$$\hat w=w\cdot\underbrace{s_{b_j-e}\cdot\dots\cdot s_{b_{j-1}+1}}
\cdot\underbrace{s_{b_j-e+1}\cdot\dots\cdot s_{b_{j-1}+2}}
\cdot\dots\cdot\underbrace{s_{b_j-1}\cdot\dots\cdot s_{b_{j-1}+e}}.\tag 6.9$$
If (6.8) holds and $\bold{D}_{\hat w}$ is irreducible, then $\rho $ 
maps $\bold{D}_w$ birationally onto $\bold{D}_{\hat w}$.

\endproclaim
\demo{Proof} Note first that if (6.8) is satisfied, then (6.9) follows directly from the definition of $\hat w$.

By the construction of $\hat w$, we have
$r_{\hat w}(b_j,p)=r_w(b_j,p)$, for all $1\leq p\leq n$. 
Since $w(b_j)<w(b_j+1)$ by assumption, it follows from
\cite{F1, Proposition 4.2} that by adding the conditions
$$ \operatorname{rank}(V_p\otimes\Cal O\rightarrow(\rho^*B_j)^*)\leq r_w(b_j,p),\ 1\leq p\leq n$$ 
to (6.5) and (6.6), we obtain the {\it same} locus $\bold{D}_w$ on $G_e(B_j/B_{j-1})$. In other words, $\bold{D}_w$ is contained in $\rho^{-1}(\bold{D}_{\hat w})$. Consider the Grassmann bundle obtained by
restriction
$$\rho:\rho^{-1}(\bold{D}_{\hat w})\rightarrow\bold{D}_{\hat w}.$$
It is not hard to see that $\bold{D}_w$ is a Schubert variety in this bundle,
of positive relative dimension, unless (6.8) holds, in which case it intersects
each fibre in a point, and the Lemma follows.
\quad\qed\enddemo

We will now prove the following equivalent reformulation of Theorem 3.1.

\proclaim{Theorem 3.1'}
The GW-number $\langle \Omega_{\alpha_{i,j}},\Omega_w,\Omega_v\rangle _{\overline d}$ vanishes, unless $\overline d$ is one of the multiindices
$\overline e_{\bold{h}\bold{l}}$ of Lemma 5.8, such that
$\ell (w\cdot\gamma_{\bold{h}\bold{l}})=\ell (w)-\sum_{c=1}^m(n_{l_c+1}-n_{h_c})$,
and $v$ is dual to one of the permutations $w''\cdot\delta_{\bold{h}\bold{l}}$, in which cases it is equal to 1. 
\endproclaim

\demo{ Proof of Theorem 3.1'}
The idea is to degenerate the corresponding intersection of generalized Schubert varieties, as in the previous subsection.

Specifically, let $\overline d$ be any multiindex not identically 0, 
and let $v\in S$ be such
that $c:=i+\ell (w)+\ell (v)=\operatorname{dim}H_{\overline d}$. 
Let $y,t\in\Bbb P^1$ be distinct points. 

We now claim that the intersection
$$Z:=\overline\Omega_{\alpha_{i,j}}(y)\bigcap
\overline\Omega_w(y)\bigcap\overline\Omega_v(t)\tag 6.10 $$ is either empty, or 
purely 0-dimensional. 

Indeed, by Theorem 4.4 $(i)$, it is enough to show that the restriction of $Z$ to $\Cal B_{\overline d}$ is either empty, or purely
of dimension 0. As in the proof of Theorem 4.4 $(ii)$, this reduces to showing
that the codimension in $\Cal U_{\overline e}$ of
$$\left (\pi ^{-1}(\Bbb P^1\times\overline\Omega_{\alpha_{i,j}}(y))\bigcup
\tilde\Omega_{\alpha_{i,j}}^{\overline e}(y) \right )\bigcap\left (
\pi ^{-1}(\Bbb P^1\times\overline\Omega_w(y))\bigcup
\tilde\Omega_w^{\overline e}(y) \right )\bigcap\left (
\pi ^{-1}(\Bbb P^1\times\overline\Omega_v(t))\bigcup
\tilde\Omega_v^{\overline e}(t) \right )$$ is at least
$c+1-\sum_{i=1}^ke_i(n_{i+1}-n_i)-\sum_{i=1}^ke_i(e_i-e_{i-1})$, for all multiindices $\overline e\neq (0,\dots ,0)$, satisfying (5.1) and (5.2). 
We have seen already in the proof of Theorem 4.4 $(ii)$ that
the only intersection which may be nonempty is
$$\tilde\Omega_{\alpha_{i,j}}^{\overline e}(y)\bigcap
\tilde\Omega_w^{\overline e}(y)\bigcap
\pi ^{-1}(\{ y\}\times\overline\Omega_v(t)),\tag 6.11$$ which lies inside
$\Cal U_{\overline e}(y)$. But we can rewrite (6.11) as 
$$\psi_ {\overline e}(y) ^{-1}(\Omega_{\tilde\alpha_{i,j}^{\overline e}})\bigcap
\psi_{\overline e}(y) ^{-1}(\Omega_{\tilde w^{\overline e}})\bigcap
\pi ^{-1}(\{ y\}\times\overline\Omega_v(t)),\tag 6.12$$
where $\psi_ {\overline e}(y):\Cal U_{\overline e}(y)\rightarrow\tilde F
_{\overline e}$ is the morphism of Remark 5.6 $(ii)$.
The codimension of (6.12) in $\Cal U_{\overline e}(y)$
satisfies the required estimate by  
Kleiman's transversality Theorem, Lemma 5.7, and Lemma 5.8
$(i)$.

The conclusion of Proposition 6.3 applies therefore in this case, and the GW-invariant $\langle \Omega_{\alpha_{i,j}},\Omega_w,\Omega_v\rangle _{\overline d}$ can be computed
as the degree
of $[Z]$ in the Chow ring of the hyperquot scheme. 
But we know even more! Namely, if $Z$ is nonempty,
all the inequalities we have used above must be in fact equalities. By
Lemma 5.8 $(ii)$ and $(iii)$, this implies that $Z$ is contained in
the (disjoint!) union of "strata"
$$\bigcup_{\overline e_{\bold{h}\bold{l}}}\phi _{\overline e_{\bold{h}\bold{l}}}
\left (\Cal U_{\overline e_{\bold{h}\bold{l}}}(y)\right ),$$
where the union is over all $\overline e_{\bold{h}\bold{l}}$, such that
$\ell (w\cdot\gamma_{\bold{h}\bold{l}})=\ell (w)-\sum_{c=1}^m(n_{l_c+1}-n_{h_c})$, and  for each 
$\overline e_{\bold{h}\bold{l}}$ as above the preimage
$\phi _{\overline e_{\bold{h}\bold{l}}}^{-1}(Z)$ is given
by the intersection (6.12), with $\overline e$ replaced by
$\overline e_{\bold{h}\bold{l}}$. 

At this point we need the following
\proclaim{Lemma 6.5} The intersection 
$$\psi_{\overline e_{\bold{h}\bold{l}}}(y)^{-1}(\Omega_{\tilde\alpha_{i,j}^{\overline e_{\bold{h}\bold{l}}}})
\bigcap\psi_{\overline e_{\bold{h}\bold{l}}}(y)^{-1}(\Omega_
{\tilde w^{\overline e_{\bold{h}\bold{l}}}})\bigcap
\pi ^{-1}(\{ y\}\times\overline\Omega_v(t))$$
is empty whenever
$\overline d\neq\overline e_{\bold{h}\bold{l}}$. \endproclaim

Granting this for a moment, let's complete the proof of Quantum Pieri.

Recall that $\Cal U_{\overline e_{\bold{h}\bold{l}}}(y)$ can be realized as a succesion of Grassmann bundles over an open subscheme 
$\Cal V\subset\{ y\}\times
\Cal H\Cal Q_{\overline d-\overline e_{\bold{h}\bold{l}}}$ (cf. the proof of Lemma 5.1). The above claim says that $Z$ is empty, except possibly when
$\overline d$ is one of the multiindices $\overline e_{\bold{h}\bold{l}}$
described above. In this case, 
$$\{ y\}\times\Cal H\Cal Q_{\overline d-\overline e_{\bold{h}\bold{l}}}=\Cal V=
\{ y\}\times H_{\overline d-\overline e_{\bold{h}\bold{l}}}=\{ y\}\times F,$$
and $\Cal U_{\overline e_{\bold{h}\bold{l}}}(y)$ is {\it projective}. 
Moreover, the map $\phi_{\overline e_{\bold{h}\bold{l}}}(y):
\Cal U_{\overline e_{\bold{h}\bold{l}}}(y)\rightarrow
\Cal H\Cal Q_{\overline d}$ is an embedding, by Theorem 5.2 $(ii)$. 
It follows that the degree of $[Z]$
in the Chow ring of $\Cal H\Cal Q_{\overline d}$ is given by
$$\int_{\Cal U_{\overline e_{\bold{h}\bold{l}}}(y)}\psi_{\overline e_{\bold{h}\bold{l}}}(y)^*[\Omega_{\tilde\alpha_{i,j}^{\overline e_{\bold{h}\bold{l}}}}]\cup\psi_{\overline e_{\bold{h}\bold{l}}}(y)^*[\Omega_
{\tilde w^{\overline e_{\bold{h}\bold{l}}}}]\cup\pi^*[\Omega_v].\tag 6.13$$
By applying the classical Pieri formula (Theorem 1.5) on 
$\tilde F_{\overline e_{\bold{h}\bold{l}}}$, we can rewrite (6.13) as
$$\sum_{\tilde w^{\overline e_{\bold{h}\bold{l}}}\overset{\tilde\alpha_{i,j}^
{\overline e_{\bold{h}\bold{l}}}}\to{\longrightarrow}w''}
\int_{\Cal U_{\overline e_{\bold{h}\bold{l}}}(y)}
\psi_ {\overline e_{\bold{h}\bold{l}}}(y) ^*[\Omega_{w''}]\cup\pi^*[\Omega_v].$$

The subscheme $\psi_ {\overline e_{\bold{h}\bold{l}}}(y) ^{-1}(\Omega_{w''})$
is the degeneracy locus inside $\Cal U_{\overline e_{\bold{h}\bold{l}}}(y)$
determined by
$$\operatorname {rank}\ (V_p\otimes \Cal O\rightarrow K_q^*)
\leq r_{w''}(q,p),\ p=1,\dots ,n,\ q\in\{ m_1,\dots ,m_a\} .\tag 6.14$$

By Kleiman's transversality theorem, we may assume that both
$\psi_ {\overline e_{\bold{h}\bold{l}}}(y) ^{-1}(\Omega_{w''})$ and the intersection
$\psi_ {\overline e_{\bold{h}\bold{l}}}(y) ^{-1}(\Omega_{w''})\bigcap
\pi ^{-1}(\Omega_v)$ have the expected codimension. Hence
$$\psi_ {\overline e_{\bold{h}\bold{l}}}(y) ^*[\Omega_{w''}]\cup\pi^*[\Omega_v]=
[\psi_ {\overline e_{\bold{h}\bold{l}}}(y) ^{-1}(\Omega_{w''})\bigcap
\pi ^{-1}(\Omega_v)]$$ in the Chow ring of 
$\Cal U_{\overline e_{\bold{h}\bold{l}}}(y)$.
Recall that $\pi :\Cal U_{\overline e_{\bold{h}\bold{l}}}(y)\rightarrow F$
can be realized as a succesion of Grassmann bundle projections
(cf. the proof of Lemma 5.1). 
By applying Lemma 6.4 $(i)$ to each of these Grassmann bundles,
starting from the top, we get that the image of 
$\psi_ {\overline e_{\bold{h}\bold{l}}}(y) ^{-1}(\Omega_{w''})$ under the projection $\pi $ is the Schubert variety $\Omega_{w'''}\subset F$,
where $w'''$ is the permutation (in $S$ !) obtained from $w''$ by the succesive
applications of Lemma 6.4 $(i)$. By Lemma 6.4 $(ii)$, it follows that
$\pi_*[\psi_ {\overline e_{\bold{h}\bold{l}}}(y) ^{-1}(\Omega_{w''})]=0$, 
unless the condition (6.8) is satisfied in every instance
where we have used Lemma 6.4 $(i)$, in which case 
$\pi_*[\psi_ {\overline e_{\bold{h}\bold{l}}}(y) ^{-1}(\Omega_{w''})]=
[\Omega_{w'''}]$.
 
Moreover, if this happens, the permutation
$w'''$ is obtained from $w''$ by applying succesively the receipe (6.9).
Using the fact that the simple transpositions $s_i$ and $s_j$ commute whenever
$i$ and $j$ are not consecutive integers, it follows easily that
$w'''=w''\cdot\delta_{\bold{h}\bold{l}}$ and
$$\ell (w''')=\ell (w''\cdot\delta_{\bold{h}\bold{l}})=\ell (w'')-m-\sum_{c=1}^m(n_{l_c}-n_{h_c-1}).$$
From the projection formula
$$\int_F\pi_*[\psi_ {\overline e_{\bold{h}\bold{l}}}(y)^{-1}(\Omega_{w''})
\bigcap\pi^{-1}(\Omega_v)]=
\int_F[\Omega_{w''\cdot\delta_{\bold{h}\bold{l}}}]\cup[\Omega _v].\tag 6.15$$

By Theorem 1.2 the latter intersection number vanishes, unless $v$ is the
permutation in $S$ dual to $w''\cdot\delta_{\bold{h}\bold{l}}$, in which case
it is equal to 1. This implies that the same holds for the intersection number (6.13).

Summarizing, $\operatorname{deg}[Z]$ vanishes, unless all the conditions 
stated in Theorem 3.1' are satisfied, in which case it is equal to 1, and
moreover, we have seen that $\operatorname{deg}[Z]=\langle \Omega_{\alpha_{i,j}},\Omega_w,\Omega_v\rangle _{\overline d}$. This completes
the proof of Quantum Pieri. \quad\qed\enddemo

\demo{Proof of Lemma 6.5} For simplicity, we will omit
$\bold{h}$ and $\bold{l}$ from the notation. We recall first the situation we're dealing with. There is a diagram
$$\CD \Cal U_{\overline e}(y) @>\psi_ {\overline e}(y) >>
\tilde F_{\overline e}\\ @V\pi VV @.\\ \Cal V @.\endCD $$ with
$\{ y\}\times H_{\overline d-\overline e}\subset\Cal V
\subset\{ y\}\times\Cal H\Cal Q_{\overline d-\overline e}$
and $\pi $ a composition of Grassmann bundle projections.
Let $$W:=\tilde\Omega_{\alpha_{i,j}}^{\overline e}(y)\bigcap
\tilde\Omega_w^{\overline e}(y)=\psi_{\overline e}(y)^{-1}(\Omega_{\tilde\alpha_{i,j}^{\overline e}})
\bigcap\psi_{\overline e}(y)^{-1}(\Omega_
{\tilde w^{\overline e}}),$$ 
We may assume that $W$ is irreducible, of the expected codimension
$\ell (\tilde\alpha_{i,j}^{\overline e})+\ell(\tilde w^{\overline e})$, while the intersection 
$W\bigcap\pi ^{-1}(\{ y\}\times\overline\Omega_v(t))$ is a nonempty finite set
consisting of reduced points, and supported on $\pi^{-1}(\{ y\}\times H_{\overline d-\overline e})$.
It follows then that 
$\pi (W)\bigcap(\{ y\}\times\Omega_v(t))$ is a nonempty zero-dimensional subscheme of $\{ y\}\times H_{\overline d-\overline e}$.
By Lemma 6.1, this would imply $\overline d=\overline e$,
and therefore conclude the proof,
if we can show that $\pi (W)\bigcap(\{ y\}\times H_{\overline d-\overline e})$ 
is of the form $Y(y)$, for some subvariety
$Y\subset F$. Set $Y:=ev_y(\pi (W))$, where $ev_y$ is the restriction of the
evaluation map to $\{ y\}\times H_{\overline d-\overline e}$.
Then $\pi (W)\subset Y(y)$. To get the reverse inclusion, it suffices to show 
that if there exists a map $f:\Bbb P^1\rightarrow F$ with $[f]\in \pi (W)$,
then for every $g:\Bbb P^1\rightarrow F$ such that $g(y)=f(y)$ we have 
$[g]\in \pi (W)$ as well. The map $f$ is represented by a sequence of subbundles
$$S_1\subset S_2\subset\dots\subset S_k\subset V^*\otimes\Cal O_{\Bbb P^1}.$$
By assumption, there exists a point in $W\subset\Cal U_{\overline e}(y)$, 
lying over $[f]$. This is equivalent to saying that for every $i\in\{ 1,\dots k\}$ there exist quotients $$S_i(y)\twoheadrightarrow \Bbb C^{e_i}\tag *$$ of the fibres
at $y$, together with compatible maps $\Bbb C^{e_i}\rightarrow\Bbb C^{e_{i+1}}$,
and which satisfy the degeneracy conditions defining $\tilde\Omega_{\alpha_{i,j}}^{\overline e}(y)$ and $\tilde\Omega_w^{\overline e}(y)$. If $g$ is another map and $g(y)=f(y)$, then the flag of fibres at $y$ for the sequence of subbundles corresponding to $g$ coincides with
$$S_1(y)\subset S_2(y)\subset\dots\subset S_k(y)\subset V^*\otimes\Cal O_y.$$
Hence we can take the {\it same} quotients $(*)$ to obtain a point in $W$
lying over $[g]$, and the Lemma is proved. \quad\qed\enddemo

\phantom{X}

\subhead{6.4 Proof of Orthogonality}\endsubhead

\phantom{X}

\proclaim{Lemma 6.6} Let $v_1,\dots ,v_N\in S$ and
$w_1,\dots w_M\in S$ be two collections of 
permutations satisfying the conditions

$(i)$ each $v_m$ (respectively, $w_m$) , $1\leq m\leq N$ 
(respectively, $1\leq m\leq M$)
is a cycle $\alpha_{i,j}$, for some $i$ and $j$;

$(ii)$ for each $j$, the number of cycles $\alpha_{i,j}$ among the $v_m$'s 
(respectively, $w_m$'s) is at most $n_{j+1}-n_j$;

$(iii)$ $\sum_{i=1}^N\ell (v_i)+\sum_{j=1}^M\ell (w_j)>\ell (w^{\circ})=
\operatorname{dim}F$.

Then $\langle\Omega_{v_1},\dots ,\Omega_{v_N}, 
\Omega_{w_1},\dots ,\Omega_{w_M}\rangle_{\overline d}=0$,
for every $\overline d$.
\endproclaim

\demo{ Proof} The condition $(iii)$ gives the result for $\overline d=
(0,\dots ,0)$, hence we may assume that $\overline d$ is not identically 0, and
that $\sum_{i=1}^N\ell (v_i)+\sum_{j=1}^M\ell (w_j)=\operatorname{dim}H_{\overline d}$. Let $y,t_1,\dots ,t_M\in\Bbb P^1$
be distinct points and let $\Omega_{v_1}\bigcap\dots\bigcap\Omega_{v_M}:=Y\subset F$.
By Lemma 5.7, and the conditions $(i)$ and $(ii)$, for every multiindex $\overline e$ we have the inequality
$$\sum_{m=1}^N\left (\ell (v_m)-\ell (\tilde v_m^{\overline e})\right )\leq
\sum_{j=1}^ke_j(n_{j+1}-n_j).$$
Using this, the same argument as in the proofs of Theorem 4.4 $(ii)$ and Lemma 6.5 shows that the intersection
$$\overline\Omega_{v_1}(y)\bigcap\dots\bigcap\overline\Omega_{v_N}(y)\bigcap
\overline\Omega_{w_1}(t_1)\bigcap\dots\bigcap\overline\Omega_{w_M}(t_M)$$
misses the boundary of $\Cal H\Cal Q_{\overline d}$. Therefore we can apply
Lemma 6.2 to conclude that
$$\langle\Omega_{v_1},\dots ,\Omega_{v_N}, 
\Omega_{w_1},\dots ,\Omega_{w_M}\rangle_{\overline d}=\int_{\Cal H\Cal Q
_{\overline d}}[\overline{Y(y)}]\cup [\overline\Omega_{w_1}(t_1)]\cup\dots 
\cup[\overline\Omega_{w_N}(t_N)].\tag 6.16$$

Now the same reasoning can be applied to the collection $w_1,\dots ,w_M$ to
reduce the integral in (6.16) to one involving only two subvarieties of $F$. By Lemma 6.1, all such intersection numbers
vanish whenever $\overline d\neq(0,\dots ,0)$. \quad\qed\enddemo 

\demo{Proof of Theorem 3.16} 
It follows at once from the orthogonality of the classical Giambelli polynomials and the definition of $P^q_w$ that it suffices to consider the case $\ell(w)+
\ell(v)>\ell (w^{\circ})$. This in turn follows if we show that 
$$\langle\langle G_{\Lambda_1\Lambda_2\dots\Lambda_k}
G_{\Psi_1\Psi_2\dots\Psi_k}\rangle\rangle=0\tag 6.17$$ 
for all partitions $\Lambda_1,\dots ,\Lambda_k,\Psi_1, \dots ,\Psi_k$ such that 
$$\mid\Lambda_1\mid+\dots+\mid\Lambda_k\mid +\mid\Psi_1\mid+\dots+\mid\Psi_k\mid>\ell (w^{\circ}).\tag 6.18$$ 
Recall that $G_{\Lambda_1\Lambda_2\dots\Lambda_k}$ was defined as the product
$G_{\Lambda_1}^{(1)}G_{\Lambda_2}^{(2)}\dots G_{\Lambda_k}^{(k)}$,
while each $G_{\Lambda_j}^{(j)}$ is itself a product 
of at most $n_{j+1}-n_j$ factors of type $G_i^j$, for various $i$'s (and the same for $G_{\Psi_1\Psi_2\dots\Psi_k}$).
But we already know by the
special case of the Quantum Giambelli formula (see Theorem 3.9 $(i)$) that
for every $i$ and $j$ the polynomial
$G_i^j$ represents the Schubert class $[\Omega_{\alpha_{i,j}}]$ in $QH^*(F)$!
Therefore the product $G_{\Lambda_1\Lambda_2\dots\Lambda_k}
G_{\Psi_1\Psi_2\dots\Psi_k}$ coincides with the {\it quantum} product obtained 
by replacing each $G_i^j$ by the corresponding $[\Omega_{\alpha_{i,j}}]$. By
Remark 4.8, the relation (6.17) is a consequence of Lemma 6.6. \quad\qed
\enddemo


\Refs
\widestnumber\key {KiMa}
\ref \key AS \by A. Astashkevich and V. Sadov\paper Quantum cohomology of
partial flag manifolds $F_{n_1,\dots ,n_k}$\jour  Comm. Math. Phys.  
\vol 170 \pages 503-528 \yr 1995\endref
\ref \key Beh \by K. Behrend \paper
Gromov-Witten invariants in algebraic geometry \jour preprint \yr 1996\endref
\ref \key BF \by K. Behrend and B. Fantechi \jour preprint \yr 1996
\paper The intrinsic normal cone \endref
\ref \key BM \by K. Behrend and Y. Manin \jour preprint \yr 1995
\paper Stacks of stable maps and 
Gromov-Witten invariants \endref
\ref \key BGG \by I. N. Bernstein, I. M. Gelfand, and S. I. Gelfand
\paper Schubert cells and cohomology of the space $G/P$\jour Russian Math. Surveys\yr 1973\vol 28\pages 1-26\endref
\ref  \key Be \by  A. Bertram\paper 
Quantum Schubert calculus \toappear \jour Adv. Math. \endref
\ref \key Bor \by A. Borel \paper Sur la cohomologie des espaces fibr\'es principaux et des espaces homog\`enes des groupes de Lie compacte
\jour Ann. of Math. (2)\vol 57 \yr 1953
\pages 115-207\endref
\ref \key C-F1 \by I. Ciocan-Fontanine \paper Quantum cohomology of flag
varieties \jour Internat. Math. Res. Notices, no. 6 \yr 1995 \pages 263-277 
\endref
\ref \key C-F2 \bysame \paper The quantum cohomology ring of flag varieties
\jour preprint\yr 1996\endref
\ref \key D \by M. Demazure \paper D\'esingularization des vari\'et\'es de
Schubert g\'en\'eralis\'ee\jour Ann. Scient. \'Ecole Normale Sup. \yr 1974
\vol 7\pages 53-88\endref
\ref\key E \by C. Ehresmann \paper Sur la topologie des certaines espaces 
homog\`enes\vol 35\yr 1934\pages 396-443\jour Ann. of Math.\endref
\ref \key FGP \by S. Fomin, S. Gelfand, and A. Postnikov \jour preprint \yr
1996\paper Quantum Schubert 
polynomials \endref
\ref \key F1 \by W. Fulton \pages 381-420 
\paper Flags, Schubert polynomials, degeneracy loci and determinantal formulas 
\yr1991 \vol 65 \jour Duke Math. Journal \endref
\ref \key F2 \bysame \book Intersection Theory\yr 1984\publ Springer Verlag
\endref
\ref\key FP \by W. Fulton and R. Pandharipande \paper Notes on stable maps
and quantum cohomology\jour Institut Mittag-Leffler Report No. 4 \yr 1996/97 \endref
\ref \key GK \by A. Givental and B. Kim \paper 
Quantum cohomology of flag manifolds and Toda lattices \yr1995 \jour Comm. Math.
Phys. \vol 168\pages 609-641\endref
\ref \key Kim1 \by B.Kim \paper Quantum cohomology of partial flag manifolds
and a residue formula for their intersection pairing
\jour Internat. Math. Res. Notices, no. 1 \yr 1995 \pages 1-16 
\endref
\ref \key Kim2 \bysame \paper On equivariant quantum cohomology
\jour Internat. Math. Res. Notices, no. 17 \yr 1996 \pages 1-11 
\endref
\ref \key Kim3 \bysame \jour UC Berkeley Thesis \yr 1996 \paper
Gromov-Witten invariants for flag manifolds \endref
\ref \key KiMa \by A.N. Kirillov and T. Maeno \paper Quantum
double Schubert polynomials, quantum Schubert polynomials and 
Vafa-Intriligator formula \jour preprint \yr 1996 \endref
\ref \key Kl \by S. Kleiman \pages 287-297 \paper The transversality of a 
general translate \yr 1974 \vol 38 \jour Compositio Math. \endref
\ref \key Kon \by M. Kontsevich \paper Enumeration of rational curves via
torus actions \inbook in {\it The moduli space of curves}, R. Dijkgraaf,
C. Faber and G. van der Geer, eds. \publ Birkhauser \yr 1995
\pages 335-368\endref
\ref \key KM \by M. Kontsevich and Y. Manin \paper Gromov-Witten classes,
quantum cohomology and enumerative geometry \yr 1994 \jour Comm. Math. Phys.  
\vol 164 \pages 525-562 \endref
\ref \key LS1 \by A. Lascoux and M.-P. 
Sch\"utzenberger \paper Polyn\^omes de Schubert \pages 447-450 \yr 1982
\vol 294 \jour C.R. Acad. Sci. Paris \endref
\ref \key LS2 \bysame \paper Symmetry and flag manifolds \inbook in 
{\it Invariant Theory},  F. Gherardelli ed.,  Lecture Notes in Math. \vol 996 
\publ Springer \publaddr Berlin \yr 1983 \pages 118-144 \endref
\ref\key Lau \by G. Laumon\paper Un anlogue global du c\^one nilpotent \vol 57
\yr 1988 \jour Duke Math. Journal \pages 647-671\endref
\ref \key LT1 \by J. Li and G. Tian\paper The quantum cohomology of homogeneus
varieties\jour J. Algebraic Geom.\toappear\endref 
\ref \key LT2 \bysame\paper Virtual moduli cycles and Gromov-Witten invariants
\yr 1996 \jour preprint \endref
\ref \key M\by I. G. Macdonald \book Notes on Schubert Polynomials 
\publ LCIM, D\'epartement de math\'ematiques et d'informatique, 
Universit\'e du Qu\'ebec \`a Montr\'eal \yr 1991 \endref 
\ref \key MS \by D. McDuff and D. Salamon \book J-holomorphic curves
and quantum cohomology \publ Univ. Lecture Ser. 6, Amer. Math. Soc.
\publaddr Providence, RI\yr 1994\endref
\ref \key Pe \by D. Petersen \jour unpublished \yr 1996\endref
\ref \key RT \by Y. Ruan, G. Tian \paper 
A mathematical theory of quantum cohomology \pages 259-367 \yr 1996 \vol 42 
\jour J. Differential Geom. \endref
\ref \key ST \by B. Siebert and G. Tian \paper 
On quantum cohomology of Fano manifolds and a formula of Vafa and Intriligator
\yr1994 \jour preprint \endref
\ref \key S \by F. Sottile \paper Pieri's rule for flag manifolds and 
Schubert polynomials\jour Annales de L'Institut Fourier\vol 46\yr 1996
\pages 89-110\endref
\ref\key T\by G. Tian\paper Quantum cohomology and its associativity
\jour preprint\yr 1995\endref
\ref \key V \by C. Vafa \paper Topological 
mirrors and quantum rings \inbook in {\it Essays on Mirror Manifolds},  
S.T. Yau ed. \publ International Press \publaddr Hong Kong \yr 1992 \endref
\ref \key W \by E. Witten \paper Topological
sigma model \pages 411-449 \yr 1988 \vol 118 \jour Commun. Math. Phys. \endref

 

\endRefs
\enddocument